\magnification 1200
\def\R{{\rm I\kern-0.2em R\kern0.2em \kern-0.2em}}
\def\N{{\rm I\kern-0.2em N\kern0.2em \kern-0.2em}}
\def\P{{\rm I\kern-0.2em P\kern0.2em \kern-0.2em}}
\def\B{{\rm I\kern-0.2em B\kern0.2em \kern-0.2em}}
\def\Z{{\rm I\kern-0.2em Z\kern0.2em \kern-0.2em}}
\def\C{{\bf \rm C}\kern-.4em {\vrule height1.4ex width.08em depth-.04ex}\;}
\def\B{{\bf \rm B}\kern-.4em {\vrule height1.4ex width.08em depth-.04ex}\;}

\def\D{{\Delta}}
\def\DD{{\overline \Delta}}
\def\bD{{b\Delta}}
\def\z{{\zeta}}
\def\cC{{\cal C}}

\def\L{{\cal L}}
\def\cC{{\cal C}}

\def\G{{\Gamma}}

\def\Tt{{(T-t)(T-1/t)}}
\font\ninerm=cmr8
\
\vskip 25mm
\centerline {\bf MEROMORPHIC EXTENSIONS FROM SMALL FAMILIES OF CIRCLES}
\centerline {\bf AND HOLOMORPHIC EXTENSIONS FROM SPHERES  }
\vskip 4mm
\centerline{Josip Globevnik}
\vskip 6mm
{\noindent \ninerm  ABSTRACT \  Let $\B$ be the open 
unit ball in $\C^2$ and let 
 $a, b, c$  be three points in $\C^2$ which do not lie in a complex line, 
 such that the complex line 
 through $a, b$ meets $\B$ and such that if one of the points $a, b$ is in 
 $\B$ and the other in $\C^2\setminus \overline\B$ then $<a|b>\not= 1$ and such that at 
 least one of the numbers
 $<a|c>,\ <b|c>$ is different from $1$. We prove that 
 if
 a continuous function $f$ on $b\B$ extends holomorphically into $\B$ 
 along each complex line which meets $\{ a, b, c\}$ then $f$ extends holomorphically through $\B$, 
 This generalizes the recent result of L.\ Baracco who proved such a result in the case when the
 points $a, b, c$ are contained in $\B$. The proof is quite different from the one 
of Baracco and uses 
the following one variable result which we also prove in the paper:  
  \ Let $\D$ be the open unit disc in $\C$. Given $\alpha\in\D$ 
let $\cC_\alpha $ be 
the family of all circles in $\D$ obtained as the images of 
circles centered at
the origin under 
an automorphism of $\D$ that maps $0$ to $\alpha$. \ Given 
$\alpha, \beta \in\D,\ \alpha\not= 
\beta$, and $n\in\N$,  a continuous function $f$ on $\DD$ extends
meromorphically from every circle 
$\Gamma\in\cC_\alpha\cup\cC_\beta$ through the disc bounded 
by $\Gamma $ with
the only pole at the center of $\Gamma$ of degree 
not exceeding $n$ if and only if $f$ is of the form \ $
f(z) = a_0(z)+a_1(z)\overline z +\cdots +a_n(z)\overline z^n\ \ (z\in\D)
$
\ where the functions $a_j,\ 0\leq j\leq n$, are holomorphic on $\D$.}

\vskip 6mm

\bf 1. The main results \rm
\vskip 2mm
Denote by $\D$ the open unit disc in $\C$. If $a\in\C$ and $r>0$ write $\D (a,r)
= \{ \z\in \C\colon \ |\z - a
|<r\} $. Given $\alpha\in \D$ the Moebius map
$$
\z\mapsto M_\alpha (\z )={{\alpha-\z}\over{1-\overline\alpha \z}}\ \ \ (\z\in \D )
$$
maps the circles in $\D$ centered at the origin to the circles
$$
\Bigl\{ {{\alpha-R\z}\over{1-\overline\alpha R\z}}\colon \ \z\in\bD\Bigr\},\ 0<R<1
$$
which are called the circles with the hyperbolic center $\alpha$ and we denote
this family of circles by $\cC_\alpha$. In particular, $\cC_0$ is the family of all 
circles in $\D$ centered at the origin.

If $\G$ is a circle we denote by $c(\G )$ its center, moreover, if $\G \subset\D $
we denote by $h(\G )$ its hyperbolic center, that is, the unique point such that $\G \in \cC _{h(\G )}$. In the case 
when $\alpha \in \bD$ then we denote by $\cC_\alpha $ 
the family of all circles in $\overline\D$ which pass through $\alpha $.

We say that a continuous function $f$ on a circle $\G$ extends holomorphically
(meromorphically) from $\G$ if 
it extends holomorphically (meromorphically) through the disc bounded by $\G$. 

If a function $f$ on $\D$ has the form
$$
f(z) = a_0(z)+a_1(z)\overline z+\cdots +a_n(z){\overline z}^n\ \ (z\in \D)
\eqno (1.1)
$$
where $a_i,\ 0\leq i\leq n$, are holomorphic functions on $\D$ then we say 
that $f$ is a
\it polyanalytic function \rm on $\D$ (of order $n\leq 1$ if $a_n\not\equiv 0)$. 
Polyanalytic 
functions of order zero are holomorphic functions.

If $f$ is a polyanalytic function of order $\leq n$ on $\D$ then for each circle
$\G\subset\D$, the function 
$z\mapsto (z-c(\G ))^nf(z)$ extends holomorphically from $\G$; in other words,
$f$ has a meromorphic extension 
from $\G$\ (to the disc bounded by $\G$) with 
the only pole at the center of $\G $, which is of degree $\leq n$. Indeed, 
if $z\in \bD (a,r)$ then 
$\overline z = \overline a + r^2/(z-a)$
so
$$(z-a)^nf(z)= (z-a)^n\Bigl[a_0(z)+a_1(z)\Bigl[\overline a+{{r^2}\over{z-a}}\Bigr]+
\cdots + a_n(z)\Bigl[\overline a+{{r^2}\over{z-a}}\Bigr]^n\Bigr]
$$ provides the necessary holomorphic extension through $\D (a,r)$.

We begin with a one-variable result.
\vskip 2mm
\noindent\bf THEOREM 1.1\ \it Let $\alpha, \beta\in \overline \D,\ \alpha\not= \beta,$  
and let $n\in \N\cup\{ 0\}$. Let $f$ be a continuous function on $\overline\D$. Assume that for every circle 
$\Gamma \in \cC_\alpha\cup\cC_\beta $ the function $z\mapsto (z-c(\G ))^nf(z)$ 
extends holomorphically from $\G $. Then the 
function $f$ is polyanalytic of order $\leq n$ on $\D$ and consequently for every circle $\Gamma \subset \overline\D$ the 
function $z\mapsto (z-c(\G ))^nf(z) $ extends holomorphically from $\G $.
\vskip 2mm
\noindent\rm In [A2] M.\ Agranovsky obtained, in a real analytic case, a characterization
of polyanalytic functions in terms of meromorphic extendibility from various 
families of circles. In the special case when $f$ is 
real-analytic Theorem 1.1 follows from his work. 
In the case when $n=0$, Theorem 1.1 follows from the result 
of A.\ Tumanov when $\alpha, \beta \in b\D$ [T2] and reduces to the results of the 
author when $\alpha\in \D,\beta\in\bD$ [G2] and when
$\alpha,\beta\in\D$ [G3]. 

Given $a\in\C^2$ we denote by $\L (a)$ the family of all complex lines passing through $a$ Given $b\in\C^2,\ b\not= 
a$, we denote by $\Lambda (a,b) $ the complex line passing through $a$ and $b$.
We denote by $\B$ the open unit ball in $\C^2$. Using Theorem 1.1 we prove 
\vskip 2mm
\noindent\bf THEOREM 1.2\ \it  Let $a, b$ be two points in $\C^2$ such that
$\Lambda (a,b)$ meets $\B$. Suppose 
that if one of the points $a, b$ is in $\B$ and the other in $\C^2\setminus\overline\B$ 
then  $<a|b> \not=1$. Assume that a continuous function $f$ on $b\B$ extends holomorphically into $\B$ along every 
complex line $L\in \L (a)\cup \L (b)$. Then for any $c\in \Lambda (a,b)\cap \B$, the function f extends holomorphically  
into $\B$ along any complex line $L\in\L (c)$. \rm
\vskip 2mm
\noindent L.\ Baracco [B2] proved recently a conjecture of M.\ Agranovsky [A1]:
\it If $a, b, c \in \B$ 
do not lie on a complex line then  $\L (a)\cup\L(b)\cup \L(c)$ is 
a test family for holomorphic extendibility for $C(b\B )$, that is, 
\it if $f\in C(bB)$ 
extends holomorphically into $\B$ along each complex line $L\in\L(a)\cup
\L(b)\cup \L (c)$ where the points $a, b, c$ are in $\B$ and do not lie on
a complex line then 
$f$ extends holomorphically through $\B$.\rm \ \ Another attempt to prove the 
conjecture was presented in [A3] by M.Agranovsky who later found that the  
proof in [A3] is incomplete.  
  Our proof of Theorem 1.2 provides a new, different proof of the result of Baracco. 
Indeed, Theorem 1.2 implies 
that $f$ extends into $\B$ along every complex line that meets $\B$ 
and so $f$ extends holomorphically through $\B$ [AV, S]. Since in our Theorem 1.2 
the points $a, b$ do not have to lie 
in $\B$ we get a more general result:
\vskip 4mm
\noindent\bf COROLLARY 1.3\it \ Let $a, b, c$ be points in $\C^2$ which do not lie on a complex line, 
such that $\Lambda (a, b)$ meets $\B$ and such that if one of the points $a, b$ lies in $\B$,
the other in $\C^2\setminus\overline \B$ then $<a|b>\not= 1$ and such that at least one of the numbers $<a|c>,\ <b|c>$ is 
different from $1$. Then $\L (a)\cup\L(b)\cup \L(c)$ is 
a test family for holomorphic extendibility for $C(b\B )$. \rm
\vskip 2mm
\bf\noindent Remark\ \rm For functions $f$ in  $\cC^\infty (b\B)$,\ \   $\L (a)\cup \L(b)$ above 
is a test 
family for holomorphic extendibility, that is, the complex lines through two points 
suffice [G3].
This is no more true for functions in $\cC^k (b\B )$ [G3]. Corollary 1.3 implies that if $N+1$ 
points in the open 
unit ball of $\C^N $ do not lie in a $(N-1)$dimensional complex plane then the complex lines
passing through these points form a test family for holomorphic extendibility 
for continuous 
functions on the unit sphere in $\C^N$. It is known that the complex lines 
through two points in the unit ball 
suffice 
for functions of class $\cC^\infty$ for any dimension $N$ [G3]. 

\vskip 4mm
\bf 2.\ Poles at the hyperbolic centers
\vskip 2mm \rm We begin by a simple. but important observation of M.\ Agranovsky [A2]:
\vskip 2mm
\noindent \bf PROPOSITION 2.1\ \it Let $\G\subset\D$ be a circle bounding the open disc $D$ which 
is not centered at $0$. 
The rational extension of $z\mapsto\varphi (z) = 1-|z|^2$ from $\G$ has one zero in $D$, a single zero at $h(\G )$, 
and one pole in $D$, a simple pole at $c(\G )$. \rm
\vskip 2mm
\noindent\bf Proof.\ \rm If $\G = \bD (a,r)$ then $z\in\G $ implies that $\overline z = \overline a + r^2/(z-a)$ so
$$
z\mapsto 1-|z|^2 = 1-z\Biggl[\overline a +{{r^2}\over{z-a}}\Biggr]
\eqno (2.1)
$$
is a rational function which has one pole $a$ in $D$ which proves the second statement.
To prove the first statement, write
$$
\G = \Biggl\{ {{\alpha-R\z}\over{1-\overline \alpha R\z}}\colon\ \z\in \bD\Biggr\}
$$
where $\alpha\in \D,\ \alpha\not= 0$, and $0<R<1$. If $\z\in b\D$ then 
$$
1-\Big\vert{{\alpha-R\z}\over{1-\overline\alpha R\z}}\Big\vert^2
= {{(1-R^2)(1-|\alpha |^2) \z }\over{(1-\overline\alpha R\z)(\z -\alpha R)}}
$$
which has one zero at $\z =0$ and one zero at $\z=\infty$ which implies that 
the rational extension of $\varphi$ from $\Gamma$ has one zero at $z=\alpha$ and the 
other zero at $z=1/\overline\alpha$ which proves the first statement,
\vskip 2mm
\noindent\bf Remark\ \ \rm If $\G\subset\overline\D$ meets $\bD$, that is,
if $\G\in \cC_\alpha$ with $\alpha \in \bD$ then the rational
extension of $\varphi$ from $\G$ has one pole in $D$, a single pole at $c(\G )$, no zero in $D$ and double zero at $\alpha $.
To see this, assume with no loss of generality that $\alpha = 1$ and let 
$\G = \{ (1-r)+r\z\colon \ \z\in \bD\}$ where 
$0<r<1$. Then, for $\z\in\bD$, 
$$
\eqalign{1-|(1-r)+r\z|^2  
=\ \ &1-(1-r)^2-r(1-r)\z - r(1-r)/\z - r^2 \cr
=\ \ &{\z^2(r^2-r)-2\z (r^2-r) + r^2-r}\over\z \cr
=\ \ &{r(r-1)(\z-1)^2}\over\z  \cr}
$$
which has a double zero at $\z =1$. This, in particular, 
implies that if $\Gamma \in \cC_\alpha $  the function $\z\mapsto (z-\alpha)^2/(1-|z|^2)$, 
defined an $\Gamma\setminus\{\alpha\} $ extends continuously to $\Gamma $.

Given a continuous function $f$ on $\overline\D$ and $n\in\N$,  set
$$
F(z) = {{f(z)}\over{(1-|z|^2)^n}}
$$
Given $\gamma \in\overline\D$ consider the following two conditions (H) and (C)
$$
\left.
\eqalign{
&\hbox{if\ }\gamma\in\D 
\hbox{\ then\ } (z-\gamma)^n F(z)\ \hbox{extends holomorphically from
each } \G\in \cC_\gamma \cr
&\hbox{if\ }\gamma\in\bD 
\hbox{\ then\ } (z-\gamma)^{2n} F(z)\ \hbox{extends holomorphically from
each } \G\in \cC_\gamma \cr}\right\}\eqno (H)
$$
$$
(z-c(\G ))^nf(z) \hbox{\ extends holomorphically from each \ }\G\in C_\gamma .
\eqno (C)
$$
If $\alpha\in \D$ and $\G\in\cC_\alpha$ then by the preceding discussion
$(z-c(\G ))^n f(z)$ extends holomorphically from $\G $ if and only 
if $(z-\alpha )^n F(z)$ extends holomorphically from $\G$. Moreover, 
if $\alpha\in \bD$ and $\G\in\cC_\alpha$ then 
$(z-c(\G ))^n f(z)$ extends holomorphically from $\G $ if and only if $(z-\alpha)^{2n}F(z)$
extends holomorphically from $\G $. This gives
\vskip2mm
\noindent\bf LEMMA 2.2\ \it Let $f$ be a continuous function on $\overline\D$ and let $n\in\N$. Let 
$$
F(z)= {{f(z)}\over{(1-|z|^2)^n}}\ \ (z\in \D) .
$$
Then for each $\gamma\in\overline \D$, (H) and (C) are equivalent. \rm
\vskip 2mm
\noindent Note that if $\alpha\in\D$ then $(z-\alpha)^nF(z)$ extends holomorphically
from each $\Gamma\in\cC_\alpha $ 
if and only if $z^n(F\circ M_\alpha )(z)$ extends holomorphically from each $\G\in \cC_0$.

Assume for a moment that we have proved Theorem 1.1 in the following cases

(i) $\alpha\in\D,\ \alpha\not= 0,\ \beta = 0$

(ii) $\alpha\in \bD,\ \beta = 0$

(iii) $\alpha= -1,\ \beta =1$, 

\noindent that is, that we have proved that 
if (C) holds for $\gamma =\alpha,\ \gamma=\beta$ with 
$\alpha, \beta$ as in (i), (ii), (iii), then (C) holds for every $\gamma\in\DD$, or, 
equivalently, if (H) holds for for $\gamma =\alpha,\ \gamma=\beta$ with 
$\alpha, \beta$ as in (i), (ii), (iii), then (H) holds for every $\gamma\in\DD$.

Assume now that $\alpha\in\D, \beta\in\D,\beta\not=\alpha,$ and that (C) holds for
$\gamma=\alpha, \gamma=\beta$. By Lemma 2.2 (H) holds for $\gamma =\alpha, 
\gamma =\beta$. This  means that both
$$
z^n(F\circ M_\alpha)(z)\hbox{\ \ and\ \ }z^n(F\circ M_\beta )(z) 
$$
extend holomorphically from each $\G \in \cC_0$. In particular,
$$
z^n(F\circ M_\alpha)\circ (M_\alpha^{-1}\circ M_\beta)(z)
$$
extends holomorphically from each $\G\in\cC_0$ so for some $\delta \not=0$,\ \ 
$z^n(F\circ M_\alpha)(M_\delta (z))$ extends holomorphically from every
$\G\in\cC_0$.  Since the function 
$z\mapsto (1-|z|^2)^n (F\circ M_\alpha)(z) $ is also 
continuous on $\DD$, Theorem 1.1 in the case (i) above  implies that 
for each $\omega\in\D$ the function $z^n(F\circ M_\omega )(z)$ 
extends 
holomorphically from each $\G\in\cC_0$ which means that (H) holds for every
$\gamma$ and consequently (C) holds for every $\gamma $.  This shows that 
Theorem 1.1 holds for every $\alpha\in\D$ and every $\beta\in\D ,\beta\not=\alpha$.

Assume now that (C) holds for $\gamma =\alpha,\ \gamma=\beta$ where 
$\alpha\in\D,\ \beta\in\bD$. By Lemma 2.2, (H) holds for $\gamma =\alpha,
\gamma =\beta$. In particular,\ $z^n(F\circ M_\alpha )(z)$  extends
holomorphically for each $\Gamma\in\cC_0$ 
and $(z-\beta)^{2n} F(z)$ extends holomorphically from each 
$\Gamma \in \cC_\beta$. It follows that $(z-\beta)^{2n} (F\circ M_\alpha)
(M_\alpha^{-1}(z))$ extends holomorphically 
from each $\Gamma\in \cC_\beta $ which means that 
$w\mapsto (M_\alpha (w) -M_\alpha M_\alpha^{-1}\beta )^{2n}
(F\circ M_\alpha)(w)$ extends holomorphically from each 
$\G\in \cC_{M_\alpha^{-1}(\beta)}$ which is the same as to say
that $(w-M_\alpha^{-1}(\beta )^{2n}(F\circ M_\alpha)(w)$ extends holomorphically from
each $\G\in \cC_{M_\alpha^{-1}(\beta)}$, that is, $F\circ M_\alpha$ satisfies (H) 
for $\gamma = 0$ and for $\gamma = M_\alpha^{-1}(\beta )$. By 
Theorem 1.1 in the case (ii) above 
it follows that 
(H) holds for $F\circ M_\alpha $ for every $\gamma $, that is, 
$z^n (F\circ M_\alpha )(M_\omega (z))$ 
extends from $\cC_0$ for every $\omega $ so so for every $\delta \in\D$,\ \ 
$z^n(F\circ M_\delta)(z)$ extends 
holomorphically from each $\G$ in $\cC_0$, that is, (H) and 
consequently (C) holds for every $\gamma $. 

In the same way we show that if $\alpha, \beta \in \bD$ and if (C) holds for $\gamma =\alpha,  \gamma =\beta$, then, 
using Theorem 1.1 in the case (iii) above we see that (C) holds for every $\gamma\in\DD$.

All this shows that, after a rotation if necessary,  
 it is enough to prove Theorem 1.1 in the cases when
$$
\alpha=t,\ \ 0<t<1, \hbox{\ \ and\ \ }\beta =0 ,
\eqno (2.2)
$$
$$
\alpha=1, \ \ \beta =0,\eqno (2.3)
$$
$$
\alpha = -1, \ \ \beta=1 . \eqno (2.4)
$$
\noindent Suppose that we have done this. Then we have also proved
\vskip 2mm
\noindent\bf THEOREM 2.3\ \it Let $F$ be a continuous function on $\D$ such that 
$z\mapsto (1-|z|^2)^n F(z)$ is continuous on $\DD$. Let $\alpha, \beta \in\DD,
\ \alpha\not=\beta $ and assume that 
(H) holds for $\gamma = \alpha$ and $\gamma=\beta $. Then
$$
F(z) ={{f(z)}\over{(1-|z|^2)^n}}\ \ (z\in\D) 
$$
where $f$ is polyanalytic of degree $\leq n$. In particular, (H) holds for every $\gamma\in\DD$. 
\vskip 4mm
\bf 3.\ Proof of Theorem 1.2 assuming Theorem 1.1.
\vskip 2mm\rm 
In this section we deduce Theorem 1.2 from Theorem 2.3. Along the way we describe, 
for a given complex line $\Lambda $ in $\C^2$ that meets $\B$, those continuous functions on $b\B$ 
that extend holomorphically into $\B$ along each complex line $L\in \L (c),\ c\in \Lambda\cap\B$. We
already know that each such function which is of class $\cC^\infty $ necessarily extends 
holomorphically through $\B$.

So let $a,b$ be two points in $\C^2$ such that $\Lambda (a,b)$ meets $\B$, and 
suppose that if one of the points is in $\B$, the other in $\C^2\setminus\overline\B$
then $<a|b>\not=1$. Assume that 
 $f\in C(b\B)$ extends holomorphically 
 into $\B$ along every complex line $L\in\L (a)\cup \L (b).$
 
 As in [G3] We use Moebius transforms to show that it is
 enough to prove the statement of Theorem 1.2 in the special case 
 when $\Lambda (a,b)$ is the $z$-axis and in one of the following two cases:
 
 (i) $a= (t_1, 0), b= (t_2,0)$ where $0\leq t_1<t_2<\infty,\  t_2\not= 1/t_1$
 
 (ii) $a= (1,0),\ b= (-1, 0)$.
 
 As in [G3] or [A1] we now use the Fourier series decomposition and averaging 
 to reduce the problem in $\C^2$ to a series of one variable problems.
 For each $z\in\D$, write the Fourier series
 $$
 f(z,e^{i\theta}\sqrt{1-|z|^2} \sim \sum_{n=-\infty}^\infty c_n(z)\bigl(
 e^{i\theta}\sqrt{1-|z|^2}\bigr)^n \eqno (3.1)
 $$ 
 so that 
 
 $$
  f(z,w) \sim \sum_{n=-\infty}^\infty c_n(z)w^n\ \ ((z,w)\in bB,\ w\not=0). 
  \eqno (3.2)
$$
Clearly the coefficients 
 $$
  c_n(z)= \left( {1\over{\sqrt{1-|z|^2}}}\right)^n{1\over{2\pi}}
  \int_{-\pi}^{\pi}e^{-in\theta}f(z,e^{i\theta}\sqrt{1-|z|^2}) d\theta .
 $$
 are continous on $\D$ and if $n<0$ they extend continuously
 to $\DD$ with zero boundary values. 
 
 If $z_0\in\C$ and if $f$ extends holomorphically into $B$ along each 
 complex line passing through $(z_0,0) $
 then in the sum (3.2) the same holds for each term: 
 $$
 \Psi _n(z,w) = {1\over{2\pi}}\int_{-\pi}^\pi e^{-in\theta }f(z,e^{i\theta}w)d\theta  =
 w^nc_n(z) ,
 $$
 a continuous function on $b\B$. Converse also holds: If each term $w^nc_n(z)$ in (3.2) 
 extends holomorphically into $\B$ along each complex line through $(z_0,0)$ then the same holds for $f$. 
 This is so since $f$ is uniformly continuous on $b\B$ 
 and so the family of functions $e^{i\theta}\mapsto
 f(z,e^{i\theta}\sqrt {1-|z|^2}) ,\ \ z\in\overline\D$ is uniformly equicontinuous 
 on $\bD$. The proof of Fejer's theorem [H] shows
 in such a case $f(z, e^{i\theta}\sqrt {1-|z|^2})$ is the limit of Cezaro means of the Fourier series (3.1) which is uniform with 
 respect to $z\in\overline\D$.
 
 Notice that $w^nc_n(z)$ extends holomorphically into $\B$ along 
 each $\Lambda\in\L ((t,0))$ if and only if $(z-t)^nc_n(z)$ extends holomorphically 
 from each $\G\in \cC_t$ in the case when 
 $0\leq t\leq 1$ and from each $
 \G\in\cC_{1/t}$ in the case when $1<t<\infty$ [G3].
 
 We will now apply Theorem 2.3 to each $c_n$. If $n\leq 0$ then 
 $c_n$ is continuous on $\DD$ and vanishes identically on 
 $\bD$ if $n<0$. If $n\geq0$ then we know that $(1-|z|^2)^{n/2}c_n(z)$ extends 
 continuously through $\DD$ and so the same holds for $(1-|z|^2)^{n}c_n(z)$.
 
 Assume now that $f$ extends holomorphically into $\B$ along each $L\in\L (a)\cup \L (b)$ where $a, b$ are as in (i). 
 If $n<0$ this implies that $c_n$ extends holomorphically from 
 each circle belonging to $\cC_{\tau_1}\cup \cC_{\tau_2}$ where 
 $0\leq \tau_1<\tau_2\leq 1$ and since $c_n$ extends continuously to $\DD$ 
 it follows by the main results of [G2, G3] that $c_n$ is holomorphic on $\D$.
 Since it is continuous on $\DD$ and 
 vanishes identically 
 on $\bD$ it follows that it vanishes identically on $\DD$.
 
 Now, let $n>0$. Now again, we have that $(z-\tau_j)^nc_n(z)$
 extends holomorphically from each $\G\in \cC_{\tau_j},\ j=1,2$, 
 which, by 
 Theorem 2.3 implies that
 $$
 c_n(z) ={{g_n(z)}\over{(1-|z|^2)^n}}\ \ (z\in\D)
 $$
 where the function $g_n$ is polyanalytic of order $\leq n$ on $\D$ and consequently, 
 for every $\gamma\in\D,\ $\ $(z-\gamma)^nc_n(z)$ extends holomorphically  from each $\G
 \in\cC_\gamma$ which implies that for each $n>0$ and for each $\gamma\in\D$ the function
 $w^nc_n(z)$ extends into $\B$ holomorphically along each 
 $L\in\L ((\gamma, 0))$ and hence, by the preceding discussion, the same holds for $f$.
 
 Similar reasoning applies in the case (ii). This completes the proof of Theorem 1.2.
 \vskip 2mm
 One should mention that the idea of multiplying $c_n$ with $(1-|z|^2)^n$ to
 achieve the regularity at the boundary and thus shifting the poles from the 
 hyperbolic centers 
 to the centers is due to M.\ Agranovsky [A1, A2]. 
 \vskip 2mm
 \noindent\bf Remark\ \rm On $b\B$ we have $|w|^2=1-|z|^2$ so it follows that, for $n>0$,  
 
$$
c_n(z) = {{g_n(z)}\over{w^n\overline w^n}}
$$ so
$$ w^n c_n(z)= {{g_n(z)}\over{\overline w^n}} ={{g_{n0}(z)+g_{n_1}(z)\overline z +\cdots 
+g_{nn}(z)\overline z^n}\over {\overline w^{n}}} .
\eqno (3.3)
$$
 It is easy to check directly that for each $k,\ 0\leq k\leq n$, and for each complex 
 line $L\in \L ((\gamma ,0))$ where $\gamma \in\D$, different from the $z-$axis, the function 
 $$
 (z,w)\mapsto {{\overline z^k}\over{\overline w^n}} = w^n{{\overline z^k}\over{(1-|z|^2)^n}}
 $$
 extends holomorphically into $B$ along $L$. This is so since
 $$
 z\mapsto (z-\gamma)^n{{\overline z^k}\over{(1-|z|^2)^n}}
 $$
 extends holomorphically from each circle $\G\in\cC_\gamma$. 
 In fact, the holomorphic extension of $\overline z$ from $\G$ has the pole 
 of order $1$ at $c(\G)$ and the same holds for $1-|z|^2$. On 
 the other hand, the rational extension of $1-|z|^2$ from $\Gamma$ has the only zero at $\gamma$.
 
 By the preceding reasoning we proved 
 \vskip 2mm
 \noindent\bf COROLLARY 3.1\ \it Let $f\in C(b\B)$ and let 
 
$$
f(z,w) \sim \sum_{-\infty}^\infty c_n(z)w^n
$$
be its Fourier series.  The following are equivalent:

(i) $f$ extends holomorphically ito $\B$ along each complex line that meets $\D\times\{ 0\}$

(ii) for each $n\in Z$,\  $(z,w)\mapsto w^nc_n(z)$ extends holomorphically into $\B$ 
along each complex line that meets $\D\times\{ 0\}$

(iii) if $n<0$ then $ c_n(z)\equiv 0$ 

and if $n\geq 0$ then the function $(z,w)\mapsto c_n(z)w^n$, 
continuous on $bB$ has the form
$$
c_n(z)w^n = c_{n0}(z){1\over{\overline w^n}}+ 
c_{n1}(z){{\overline z}\over{\overline w^n}} +
c_{nn}(z){{\overline z^n}\over{\overline w^n}}
$$
whwre $c_{n0},c_{n1},\cdots c_{nn}$ are holomorphic functions on $\D$. \rm
\vskip 2mm
\noindent Obviously, by Theorem 1.2, for either (i) or (ii) to hold it is enough that it holds for complex lines through 
two points of $\D\times\{ 0\}$. 
\vskip 2mm\noindent\bf Remark\ \rm When we want to construct examples we must take into account that the functions $(z,w)\mapsto w^nc_n(z)$ 
must be continuous on $b\B$. For instance, the following standard example is of this sort:
$$
{{w^2}\over{\overline w}} = w^3.{1\over{1-|z|^2}} = {1\over{\overline w^3}}(1-|z|^2)^2 
= {1\over{\overline w^3}}(1-2z\overline z+z^2\overline z^2).
 $$
 \vskip 1mm
 \noindent\bf Remark\ \rm We present another example, due to M.\ Agranovsky.  Given $\alpha,\ \beta \in \D$,\ 
 $\alpha\not=\beta $, the function 
 $$
 \varphi (z) = {{(z-\alpha)(z-\beta)}\over{1-|z|^2}}\ \ (z\in\D)
 $$
 extends holomorphically from every $\Gamma\in \cC_\alpha\cup \cC_\beta$.  
 
 Now, let $1 <t_1<t_2<\infty $. The function 
 $$
 (z,w)\mapsto \left\{
 \eqalign{
 w^3{{(z-1/t_1)(z-1/t_2)}\over{1-|z|^2}} 
 = {{w^2}\over{\overline w}}(z-1/t_1)(z-1/t_2)\ \ \ \ (w\not=0)\cr
 0\ \ \ \ (w=0)\cr}
 \right.
 $$
 is continuous on $b\B$ and  extends holomorphically into $\B$ along any complex line 
 meeting $\D\times \{ 0\}$ 
 as well as along any complex line in $\L ((t_1,0))\cup \L ((t_2,0))$, yet does not 
 extend holomorphically through $\B$. Obviously, finitely 
 many points on the z axis outside $\overline\B$ are also possible.
 \vskip 4mm
 \bf 4.\ Towards the beginning of the proof of Theorem 1.1 \rm
 \vskip 2mm
 We have already seen that it is enough to prove Theorem 1.1 in the cases 
 (2.2), (2.3) and (2.4). We first look at (2.2). In this case we have to prove
 \vskip 2mm
 \noindent\bf THEOREM 4.1\ \it Suppose that $0<t<1$ and let $n\in \N\cup \{ 0\}$.
 Suppose that $f$ is a continuous function on $\DD$ such that for each circle 
 $\Gamma\in \cC_0\cup\cC_t$, the function $f|\G$ extends 
 meromorphically through the disc bounded by $\G $ with 
 the only pole of order $\leq n$ at the center of \ $\G$. 
 Then $f$ is of the form
 $$
 f(z) = a_0(z)+a_1(z)\overline z +\cdots +a_m(z)\overline
 z^m
 $$
 where $m\leq n$ and where $a_j,\ 0\leq j\leq m$, 
 are holomorphic functions on $\D$,
 \vskip 2mm
 \rm\noindent The proof of Theorem 4.1 will be our major task. 
 We will provide a detailed proof. At the
 end we will indicate how to modify the proof to treat (2.3) and (2.4)

 As in [G2, G3] we shall use semiquadrics introduced in [AG] and [G1] to transform the problem 
 into a problem about holomorphic extensions of CR functions on surfaces consisting of these semiquadrics. 
 The principal idea to apply the reasoning of H.\ Lewy and H.\ Rossi 
 in this context is due to A.\ Tumanov 
 and is described in [T2], and, for holomorphic extensions applied in [G2, G3]. 
 There, the final step was to apply the Liouville theorem to conclude that the main function is constant in 
 the second variable. The proof is considerably more complicated in the case of meromorphic extensions, where we apply 
 the generalized Liouville theorem to conclude that the function is a polynomial in 
 the second variable. In [G2, G3] we were dealing with holomorphic extensions 
 to semiquadrics 
 where the maximum principle assured continuous dependence of extensions
 on the parameter. Now we will be 
 dealing with meromorphic extensions where we do not have the maximum principle
 so we will need the following preliminary 
 \vskip 2mm\noindent
 \bf LEMMA 4.1\ \it Let $I\subset \R$ be an interval. Let $n\in\N$ and let
 $(\z ,t)\mapsto \Phi (\z ,t)$ be a continuous function on $b\D\times I$ such 
 that for each 
 $t\in I$ the function $\z\mapsto \Phi (\z ,t)$ extends holomorphically 
 through $\D\setminus 
 \{ 0\}$ and has a pole at $0$ of degree $\leq n$.
 Denote the extension by $\tilde\Phi $. Then 
 $$
 \tilde\Phi (\z ,t)= {{d_{0}(t)}\over{\z ^n}} +\cdots + {{d_{n-1}(t)}\over\z} +
 \Theta (\z ,t)\ \ (\z\in\DD\setminus\{ 0\} )
 $$ 
 where the functions $d_j$ are continuous on $I$,\ \  $\Theta$ is 
 continuous on $\DD\times I$,  and for each $t\in I$, \ 
 $\z \mapsto \Theta (\z ,t)$ is holomorphic on $\D $ .
 \vskip 2mm
 \noindent\bf Proof.\rm\
 The function $(\z ,t)\mapsto \Psi (\z ,t) = \z ^n\Phi (\z ,t)$ is
 continuous on $b\D\times I$  and for each $t,\ \z\mapsto \Psi
 (\z ,t)$ extends holomorphically through $\D$. Denote this extension by $\tilde\Psi $. 
  We have
 $$ 
 \tilde\Psi (\z ,t) = d_0(t)+d_1(t)\z +\cdots 
 +d_{n-1}(t)\z^{n-1}+\z^n\Theta (\z,t)\ \ \ (\z\in\DD)
 $$
 where for each $j,\ 0\leq j\leq n-1,$
 $$
 d_j(t) = {1\over{2\pi}}\int_{-\pi}^\pi e^{-ij\theta}\Psi (e^{i\theta},t)d\theta 
 $$ 
 is continuous on $I$ and where for each $t\in I$, the function $\z\mapsto \Theta (\z ,t) $
 is holomorphic on $\D$.  The function
 $$
 (\z ,t)\mapsto \Theta (\z, t) = {{\Psi (\z ,t)}\over{\z^n}}- {{d_0(t)}\over{\z^{n-1}}}
 - \cdots {{d_{n-1}(t)}\over\z} \ \ \ (\z\in \bD)
 $$
 is continuous on $b\D\times I$ and so by the maximum principle it follows that it is continuous
 on $\DD\times I$. This completes the proof. 
 
 \vskip 4mm
 \bf 5.\ The geometry of semiquadrics, Part 1\rm
 \vskip 2mm
  Given $a\in\C$ and $r>0$ let
 $$\eqalign{
 \Lambda (a,r) & = \{ (z,w)\in\C^2\colon\ (z-a)(w-\overline a)=r^2,\ 0<|z-a|<r\} \cr
 &=\bigl\{ \bigl( z,\overline a+{{r^2}\over{z-a}}\bigr)\colon\ 0<|z-a|<r\bigr \}\cr}
 $$
 be the semiquadric associated with the circle $b\D (a,r)$. Write
 $\Sigma = \{ (\z, \overline\z )\colon\ \z\in\C\}$. \ $\Lambda(a,r)$  is a closed complex submanifold
 of $\C^2\setminus\Sigma $ which is attached to $\Sigma $ along
 $b\Lambda (a,r)= \{ (\z,\overline \z)\colon\ \z\in\bD (a,r)\}$. The crucial property 
 of these semiquadrics which will connect our problem with a problem in CR geometry 
 is the following:
 
 \noindent \it A
 continuous function $g$ extends meromorphically from a circle 
 $\bD (a,r)$ with the pole of 
 degree $\leq n$  at $a$ if and
 only if the
 function $G$, defined on 
 $b\Lambda (a,r)$
 by
 $G(\z,\overline\z )= 
 g(\z )\ (\z\in b\D (a,r))$ 
 extends holomorphically through
  ${\Lambda (a,r)}$ and has
  a pole of degree $\leq n$ at the point at infinity.  \rm
  In fact,
 if we denote by the letter $\tilde g$
 the holomorphic extension 
 of $g$ through $\D (a,r)\setminus\{ a\} $ the function 
 $$
 \tilde G\Bigl(z,\overline a + {{r^2}\over{z-a}}\Bigr) = \tilde g(z)\ \ 
 (z\in\DD(a,r)\setminus \{ a\})
 $$
 provides the necessary holomorphic extension of $G$ through $\Lambda (a,r)$.
 In our case, we have two families of circles;
 $$
 \cC_0 =\{ b\D (0,R)\colon\ 0<R\leq 1\}
 $$
 and
 $$
 \cC_t = \{ \bD (T,\rho (T))\colon\ 0\leq T\leq 1\}
 $$
 where 
 $$             \rho (T)=\sqrt\Tt \ \ 
 $$ 
 [G3] and consequently we will deal with two families of semiquadrics: the \it first \rm family 
 \ $\{ \Lambda (0,R),$ $ 0<R\leq 1\}$ and the \it second \rm family $\{ \Lambda 
 (T,\sqrt\Tt )\colon\ 0\leq T<t\}$. Our function $F(\z,\overline\z)= f(\z )\ (\z\in\DD)$ extends holomorphicall to 
 each semiquadric of each family and these extensions 
 will then define a CR function on a hypersurface in $\C^2$ that we study and describe 
 in detail.
 
 In the \it first \rm family the semiquadrics are pairwise disjoint. Let 
 $L$ be the union of their closures in $\C\times\C$;
 $$
 L = \bigcup_{0<R\leq 1}[\Lambda (0,R)\cup b\Lambda (0,R)] .
 $$
 In $[\C\setminus\{ 0\}]\times\C$, the set $L$ is a smooth hypersurface
 with piecewise smooth boundary consisting of $\Lambda (0,1)$ and
 $\{ (\z,\overline\z)\colon\ 
 \z\in \DD\setminus \{ 0\}$. 
 
 In the \it second \rm family all semiquadrics contain the point $(t,1/t)$ but otherwise are pairwise disjoint, that is, the sets 
 $\Lambda (T,\rho (T))\setminus\{ (t,1/t)\}$ are pairwise disjoint.
 The closure of their union in $[\C\setminus\{ t\}]\times \C$, that is, 
 the set
 $$
 N = \bigcup_{0\leq T<t}[\Lambda (T,\rho (T))\cup b\Lambda (T,\rho (T))]\setminus \{ (t,1/t)\} 
 $$
 is a smooth surface in $[\C\setminus\{ t\} ]\times \C$ with
 piecewise smooth boundary consisting of $\Lambda(0,1)$ and 
 $\{ (\z, \overline\z )\colon\ \z\in \DD\setminus \{ t\}\}  $.
 \vskip 4mm
 \bf 6.\ The geometry of semiquadrics, Part 2\ \rm
 \vskip 2mm
 We now fix $z\in \D\setminus\{ 0\}$ and describe how the semiquadrics $\Lambda (0,R),\ 0<R\leq 1$, and
  $\Lambda (T,\rho (T)),\ 0\leq T<t$, intersect the complex line $\{ z\} \times \C$. Each semiquadric intersects 
  $\{ z\}\times \C$ at at most one point. For a semiquadric from the first family the intersection $(z, v(R)) \in \Lambda (0,R)$
  exists iff $R(z)<R\leq 1$ where $R(z)=|z|$. As $R$ increases from $R(z)
  $ to $1$, the point $v(R)$ moves along the segment $J_z$ connecting 
  $\overline z$ and $1/z$ from $\overline z$ to $1/z$ [G3].
  
  Assume now that $\Im z\not= 0$. Denote by $\lambda_z$ the arc on the circle
  passing through $t, 1/t$ and $\overline z $ (and consequently passing through $1/z$) which 
  joins $\overline z$ and $1/z$ and does not contain $t$ and $1/t$. Let
  $T(z)$ be such that $z\in \bD (T(z),\rho (T(z)))$. For 
  the semiquadric from the second family the intersection
  $(z,w(T))\in\Lambda (T,\rho (T))$ exists if and only if $0\leq T<T(z)$. 
  As $T$ increases from $0$ to $T(z)$, $ w(T) $ moves along 
  $\lambda_z$ from $1/z$ to $\overline z$ [G3]. 
  
  Now, let $z=\eta,\ 0<\eta<t$. Let $v(R)$ be such that $(\eta, v(R))\in\Lambda (0,R)$ and 
  let $w(T)$ be such that $(\eta, w(T))\in\Lambda (T,\rho (T))$. It is easy too see that 
  as $R$ increases from $R(\eta )=\eta $ to $1$,\ $v(R)\in\R $ increases from $\eta$ 
  to $1/\eta $.  As $T$ increases from $0$ to $\eta $,\  $w(T)\in\R$
  increases from $\eta $ to $1/\eta $. As $T$ increases from $0$ 
  to $\eta $, $w(T)$ increases from $1/\eta $ to $\infty $. As $T$ 
  increases from $\eta $ 
  to $T(\eta)$ then $w(T)\in\R $ increases from $-\infty $ to $\eta $.
  
  Situation is entirely different if $z=\eta$ where either 
  $-1<\eta<0$ or $t<\eta<1$. In this case we have "folding": 
  when $R$ increases from
  $R(\eta )$ to $1$,\ $v(R)$ moves from $\eta $ to $1/\eta
  $ and when $T$ increases from $0$ to $T(\eta )$ then $w(T)$ moves
  from $1/\eta $ to $\eta $.
  
  The preceding discussion offers a "fiber" description of $L$ and $N$. We clearly have
  $$
 L=\bigcup_{z\in\DD\setminus\{ 0\}}\{ z\} \times J_z .
 $$ 
 To describe $N$, recall that so far we have only defined $\lambda_z$ for
 $z\in\DD,\ \Im z\not=0$. We now define $\lambda _z$ for 
 $z=\eta\in\R$. If $0<\eta<t$ then let 
 $\lambda_\eta = (-\infty,\eta]\cup [1/\eta ,\infty )$ and 
 when $-1\leq \eta<0$ then let $\lambda _\eta =[1/\eta,\eta]$ 
 and when $t<\eta\leq 1$ then let $\lambda _\eta = [\eta, 1/\eta ]$,
 Finally, let $\lambda_0 = (-\infty,0]$, and $\lambda_t =\{ t\}$.  Then it is easy to see that 
 $$
 N=\bigcup _{z\in\DD\setminus \{ t\}} \{ z\} \times\lambda_z .
 $$
 \vskip 4mm
 \bf 7.\ The domains bounded by $L$ and $N$ \rm
 \vskip 2mm
 Both $L$ and $N$ are smooth, unbounded surfaces with boundary. They
 have common boundary consisting of $\{ (\z,\overline \z)\colon\ 
 \z\in \DD\setminus \{ 0, t\}\} $ and $\Lambda (0,1)$. They have
 no other common points of the form $(z,w)$ where $z\in\DD\setminus [(-1,0]\cup[t,1)]$. 
  The situation is different for the points
 $z\in (-1,0)\cup (t,1)$. In this case $J_z=\lambda_z$.
 If $\eta\in (-1,0)\cup  (t,1)$ then $J_\eta =\lambda_\eta $ is the segment 
 joinig $\eta$ and $1/\eta$. .
 
 For each $z\in\D,\ \Im z\not= 0, J_z\cup \lambda _z$  is a simple
 closed curve bounding a domain that we
 denote by $D_z$. If $\Im z >0$ then $D_z$ is contained 
 in the lower halfplane and if $\Im z <0$ then $D_z$ is 
 contained in the upper halfplane. When $z\in\D$ 
 approaches a point $a\in\bD, \ \Im a\not=0$, then the domains 
 $D_z$ shrink to the point $\{ a\} $. 
 When $z,\ \Im z >0$, approaches a point $\eta ,\  
 0<\eta<t$, then $D_z$ become larger and larger and in the limit they become
 the lower halfplane. 
 When $z,\ Imz <0$, approaches a point $\eta,\  
  0<\eta<t$, then $D_z$ become larger and larger and in the limit they become
  the upper halfplane. When $z,\ \Im z>0$, approaches  a point $\eta\in (-1,0)\cup (t,1)$ 
  then the domains $D_z$ 
  become thinner and thinner and in the limit they shrink to the segment with 
  endpoints $\eta$ and $1/\eta$. The same takes place for $\Im z <0$. 
  
  As usual, denote by $\pi $ the projection $\pi (z,w)=z$.
  
  We now define two domains in $\C^2\colon\ $ Let 
  $$
  \Omega_+ = \bigcup_{z\in\D,\ \Im z>0} \{ z\} \times D_z   \ \ \ \hbox{\ and\ }\ \ \  
    \Omega_- = \bigcup_{z\in\D,\ \Im z<0} \{ z\} \times D_z   .
  $$
  Clearly $\Omega_+ \subset \{ (z,w)\colon\  \Im w <0\} $ and
  $\Omega_- \subset \{ (z,w)\colon\  \Im w >0\} $.  Obviously
  $$
  b\Omega_+ \cap \{ (z,w)\colon\  \Im z>0\} = (N\cup L)\cap \{ (z,w)\colon\  \Im z>0\} 
  $$
  $$
  b\Omega_- \cap \{ (z,w)\colon\  \Im z<0\} = (N\cup L)\cap \{ (z,w)\colon\  \Im z<0\} . 
  $$
  Further,
  $$
  b\Omega_+\cap \pi^{-1}((0,t)) = \bigcup_{0<z<t} [\{ z\} \times \{ w\in\C\colon 
  \Im w\leq 0\} ]
  $$
  and
  $$
    b\Omega_-\cap \pi^{-1}((0,t)) = \bigcup_{0<z<t} [\{ z\} \times \{ w\in\C\colon 
    \Im w\geq 0\} ].
  $$
  Obviously, 
  $$
  b\Omega_+\cap\pi^{-1}((-1, 0)) = b\Omega_-\cap\pi^{-1}((-1, 0)) =
  \bigcup_{-1<\eta<0}\{ \eta\} \times [1/\eta, \eta]
  $$
  and
  $$
  b\Omega_+\cap\pi^{-1}((t, 1)) = b\Omega_-\cap\pi^{-1}((t, 1)) =
  \bigcup_{t<\eta<1}\{ \eta\} \times [\eta, 1/\eta ] .
  $$
  \vskip 4mm
  \bf 8.\ Our CR function on $L\cup N$ \rm
  \vskip 2mm
  Let $n\in\N\cup\{ 0\}, \ 0<t<1$ and suppose that $f$ is a continuous function on
  $\DD$ that satisfies the assumptions of Theorem 4.1.
  
  Define 
  $$
  F(\z ,\overline\z )=f(\z )\ \ (\z\in\DD).
  $$
  $F$ is a continuous function on $\{ (\z ,\overline \z)\colon\ \z\in\DD\} $. The assumptions together with the continuity 
  of $F$ and Lemma 4.1 imply that the function $F$ extends from 
  $L\cap \Sigma = \{ (\z ,\overline \z)\colon \z\in\DD
  \setminus\{ 0\}\} $ continuously to $L$ such that the extension 
  $\Phi $ is holomorphic on each fiber $\Lambda (0,R)$ of $L$ with a pole of degree 
  $\leq n$ at infinity. Thus, $\Phi $ is a CR function. Similarly, the function $F$ extends from $N\cap\Sigma = \{ (\z,\overline\z)\colon\ \z\in\DD
  \setminus\{ t\} \}$ continuously to $N$ so that the extension $\Psi $ is holomorphic on each fiber
  $\Lambda (T,\rho (T))\setminus \{ (t,1/t)\},\ 0\leq T<t$, with a pole of degree $\leq n$ at 
  infinity. Thus, $\Psi $ is a CR function. The functions $\Phi $ and $\Psi $ coincide on 
  $\{(\z,\overline \z)\colon\ \z\in\DD\setminus\{  0,t\} \} \cup \Lambda(0,1)$, the common boundary of $L$ and $N$. Since 
  in $\pi ^{-1}\bigl( \D\setminus[(-1,0]\cup [t,1)]\bigr)$ 
  there are no other common points  of $L$ and $N$ we can define 
  the function $F$ on $M$ where 
  $$
  M=(L\cup N)\setminus \left[ [(-1,0]\cup [t,1)]\times \C\right]
  $$
  which extends the original $F$ and is equal to $\Phi $ on 
  $L$ and to $\Psi $ on $N$. The function $F$ so obtained is continuous on $M$ 
  and holomorphic on each fiber of $M$ so it is a CR function on $M$. 
  
  Our final goal will be to show that on $M$ we have
  $$
  F(z,w)=a_0(z) + a_1(z) w+\cdots +a_n(z)w^n
  $$
  where $a_j, \ 0\leq j\leq n$, are holomorphic functions on $\D$ which will then imply that
  $$
  f(z) = F(z,\overline z) = a_0(z) + a_1(z) \overline z +\cdots +a_n(z){\overline z}^n
  $$
  on $\D $, which is the statement of Theorem 4.1.
  
  In exactly the same way as in [G2, G3], following an idea of A.Tumanov [T2], 
  we now 
  use an argument of H.\ Lewy [L], extended by H.\ Rossi [R], 
  to show that 
  the function $F$ that is CR on $M$ extends
  holomorphically through $\{ z\} \times D_z$ for 
  each $z\in\D,
  \ \Im z\not= 0$ and that the extension so obtained 
  is holomorphic in $z$. In this way we obtain a function $F$ that is holomorphic
  on $\Omega_+$ 
  and on $\Omega_-$ and which extends continuously to 
  $\bigcup_{z\in\DD, \Im z >0} \{ z\} \times bD_z $, 
  a part of boundary of $\Omega_+$, and to 
  $\bigcup_{z\in\DD, \Im z <0} \{ z\} \times bD_z  $, 
  a part of boundary of $\Omega_-$.
  
  Before we proceed, notice that our function $F$ is well defined on $(0,t)\times\R$. 
  
  We will now show that the continuity of $F$ on $M$ implies that 
  $$
  F \hbox{\ extends continuously 
  to\ } \Omega_+\cup [(0,t)\times \{ w\colon \Im w\leq 0\} ]
  \eqno (8.1)
  $$
  and 
  $$
  F \hbox{\ extends continuously to\ } \Omega_-\cup [(0,t)\times 
  \{ w\colon \Im w\geq 0\} ].
  \eqno (8.2)
  $$ 
  We will also show that our assumptions imply that for each $\eta$, $0<\eta<t$
  there is a 
  constant $c(\eta )$ such that
    $$
    |F(\eta ,\z )|\leq c(\eta)(1+|\z |)^n\ \ (\z\in\C ).
    \eqno (8.3)
  $$  
  Suppose for a moment that we have done this. 
  
  Since $F$ is holomorphic on $\Omega_+\cup\Omega _-$ it 
  follows that the extended function $F$ is holomorphic on 
  $[\{ \eta\} \times\{ \Im\z>0\}]\cup [\{ \eta\} \times\{ \Im\z<0\}]$ and 
  continuous on $\{ \eta\}\times\C$ for each $\eta,\ 0<\eta<t$, so for each 
  such $\eta$, the function $F$ is holomorphic on $\{\eta\}\times\C$. 
  By (8.3) it follows that for each $\eta,\ 0<\eta<t,\ 
  w\mapsto F(\eta,w)$ is a polynomial of degree $\leq n$.
  \vskip 4mm
  \bf 9.\ Polynomial in the second variable \rm
  \vskip 2mm
  The fact that for each $\eta, \ 0<\eta<t,$ the function $w\mapsto F(\eta ,w)$ is a 
  polynomial of degree $\leq n$ implies, for instance, that for each such 
  $\eta$, 
 $$
 {{\partial^kF}
 \over{\partial w^k}}(\eta,w)\equiv 0\ \ (k\geq n+1,\ \Im w<0).
 $$
  For each $\eta,\ 0<\eta<t$, and for each $w_0,\ \Im w_0<0$, there is a 
  neighbourhood 
  $U\subset\C$ of $\eta $ 
  and a neighbourhood   $W$ of $w_0$, such that for each $k$,\ ${{
  \partial^kF}\over{\partial w^k}}$ is continuous on $\{ z\in U,\Im z\geq 0\}\times W$ (a 
  consequence of expressing the derivatives with the Cauchy integral formula) and  
  holomorphic on $\{ z\in U\colon \Im z>0\} \times W$. By the preceding discussion, for each 
  $k\geq n+1$,\ ${{\partial^kF}\over{\partial w^k}}$ vanishes identically 
  on $\{ z\in U,\ \Im z =0\}\times W$, so it follows 
  that it vanishes identically on $\{ z\in U, \Im z > 0\} \times W$, an open subset 
  of $\Omega_+$, and
  since it is holomorphic on $\Omega_+$, it follows that it vanishes identically 
  on $\Omega _+$ as $\Omega_+$ is connected.  Thus,
  $$
  {{\partial^kF}\over{\partial w^k}}(z,w) = 0\ \ (k\geq n+1)\ \ \hbox{\ for all\ }
  (z,w)\in\Omega_+ . 
  \eqno (9.1)
  $$
  In the same way we prove that (9.1) holds for all 
  $z\in \Omega_-$. Recall that for each $z\in\D,\ \Im z\not=0$, 
  $D_z$ is connected and hence (9.1) implies that for each such $z$, there are numbers 
  $a_0(z),\cdots ,a_n(z)$ such that 
  $$
  F(z,w) = a_0(z)+a_1(z)w+\cdots +a_n(z)w^n\ \ (w\in D_z).
  $$
  Since $F$ is holomorphic on $\Omega_+$ it follows 
  that $z\mapsto a_j(z)\ (0\leq j\leq n)$ are holomorphic on $\{ z\in\D,\ \Im z>0\}$. In the same 
  way we see that $a_0,\cdots ,a_n$ are holomorphic on $\{ z\in\D\colon \Im z <0\}$.
  
  Choose distinct points $w_0, w_1,\cdots w_n\in  \{\Im w <0\}$. 
  We have assumed (8.1)  hence $F$ extends continuously to $\Omega_+\cup
  [(0,t)\times\{ \Im w\leq 0\}]$ and it follows that
  there is a neighbourhood $U$ of $(0,t)$ in $\C$ such that 
  for each $j,\ 0\leq j\leq n $
  , the function 
  $z\mapsto F(z,w_j)$ extends continuously from $\{ z\in U,\ \Im z>0\}$ to 
  $\{ z\in U,\ \Im z>0\} \cup (0,t)$.
   Since
  $$
  \left[ \eqalign{a_0(z)\cr \cdots\cr a_n(z)\cr
  }\right] =
  \left[\eqalign{1, w_0,\cdots , w_0^n\cr \cdots\cr 
  1, w_n ,\cdots , w_n^n\cr}\right]^{-1}
  .
  \left[\eqalign{F(z,w_0)\cr \cdots \cr F(z,w_n)\cr}\right]\eqno (9.2)
  $$
  where the inverse of the Vandermonde matrix on the right exists since $w_j,\ 
  0\leq j\leq n$ are distinct, it follows that each $a_j$ extends continuously from 
  $\{ z\in\D\colon\ \Im z>0\} $ to $\{ z\in\D\colon \Im z>0\}
  \cup (0,t)$. In the same way 
  we see that the functions  $a_0,\cdots, a_n$ 
  extend continuously from $\{ z\in\D ,\Im z<0\} $ to $\{ z\in\D ,\Im z<0\}\cup (0,t)$.
  
  Recall that for each $\eta\in (0,t)$ and $s\in\R$,\ $F(\eta ,s)$ is well defined and for 
  all $s\in\R$ we have
  $$
  F(\eta ,s)= a_0^+(\eta)+a_1^+(\eta)s+\cdots + a_n^+(\eta )s^n\ \ (s\in\R)
  \eqno (9.3)
  $$
  where $a_j^+$ are continuous extensions of $a_j$ from  $\{ \z\in\D,\ \Im \z >0\} $ and 
  $$
  F(\eta, s)=  a_0^-(\eta)+a_1^-(\eta)s+\cdots + a_n^-(\eta )s^n\ \ (s\in\R)
  \eqno (9.4)
  $$
  where $a_j^-$ are continuous extensions of $a_j$ from $\{ \z\in\D,\ \Im \z <0\} $. Now,
  (9.4) and (9.3) imply that $a_j^-\equiv a_j^- $ on $(0,t)$ for each $j,\ 0\leq j\leq n$, 
  which implies that 
  there are holomorphic functions $a_j,\ 0\leq j\leq n,$ on $\D\setminus [ (-1,0]\cup [t,1)]$,
  such that 
  $$
  F(z,w)=a_0(z)+a_1(z)w+ \cdots +a_n(z)w^n\hbox{\ \ for all\ \ } 
  (z,w)\in\Omega_-\cup\Omega_+\cup [(0,t)\times \C].
  $$
  
  In the next section we show first that the functions 
  $a_j,\ 0\leq j\leq n,$ extend 
  holomorphically also across the intervals $(-1,0)$ and $(t,1)$ and then we show that the 
  remaining points $0$ and $t$ are removable 
  singularities for all these functions.
  \vskip 4mm
  \bf 10.\ Analyticity of the coefficients\rm
  \vskip 2mm
  Our extended function $F$ is well defined on both
  $$
  \bigcup_{z\in\D_+}\{ z\}\times \overline D_z\hbox{\ and\ }
  \bigcup_{z\in\D_-}\{ z\}\times \overline D_z
  $$ where $\D_+=\{ \z\in\D\colon\ \Im z>0\} $ is the upper half of $\D$ and $
  \D_-$ is the lower half of $\D$.
  
  For each $z\in\D,\ \Im z\not=0,$ denote by $I_z$ the segment from $\overline z$ 
  to $1/z$ (the straight part of $bD_z$) 
  and by $J_z$ the circular arc part of $bD_z$. The sets $I_z$ and $J_z$ 
  meet at the points $\overline z$ and $1/z$. For $-1<x<0$
  we define $I_x=J_x= (1/x, x)$. Since the points of $\{ x\}\times I_x$ belong to both families of semiquadrics there is 
  a problem of defining $F$ on $\{ x\}\times I_x$. At each point of $\{ x\} \times I_x$ we have two different 
  values of $F$: one comes from values of extensions on the first family and the other from 
  the values of extensions on the second family. From the definition of $F$ we know that
  $$
  \hbox{the function\ } F\hbox{\ extends continuously to\ }
  \bigcup_{z\in\D_+\cup\D _-\cup (-1,0)} \{ z\} \times I_z .
  \eqno (10.1)
  $$
  For $z=x\in (-1,0)$ denote this extension by $\varphi (x,w)\ \ (w\in (1/x,x))$. \ 
   Clearly $\varphi (x,w)$ is the value of the extension of the original $F$ to the semiquadric $\Lambda (0,R)$ 
   of the first family which passes through $(x,w)$, at the point $(x,w)$.
   Similarly, 
   $$
     \hbox{the function\ } F\hbox{\ extends continuously to\ }
     \bigcup_{ z\in\D_+\cup\D _-\cup (-1,0)}  \{ z\} \times J_z.
     \eqno (10.2)
  $$
  For $z=x\in (-1,0)$ denote the extension by $\psi (x,w)\ (w\in (1/x,x)$.
     Here $\psi (x,w)$ is the value of the extension of original $F$ to 
     the semiquadric $\Lambda (T,\rho (T))$ 
   of the second family which passes through $(x,w)$, at the point $(x,w)$.
   
   We will show that $\varphi$ and $\psi$ coincide and that for 
   each $x\in (-1,0)$ there are numbers $p_0(x),\cdots ,p_n(x)$, such that 
   $$
   \varphi (x,y)=\psi (x,y) =p_0(x)+p_1(x)y\cdots +p_n(x)y^n\ \ (1/x<y<x),
   $$
   and we will then show that for each $x\in (0,1)$ and for each $
   j,\ 0\leq j\leq n$,  
   $p_j(x)$ is the value at $x$ 
   of the continuous extension of $a_j$ from $\D_+\cup\D_-$ to $\D_+\cup\D_-\cup\{ x\}$.
   
   We need the following
   \vskip 2mm
   \noindent\bf LEMMA 10.1\ \it Let $p_m(w) = a_{m0}+a_{m1}w+\cdots + a_{mn}w^n$ be a sequence of polynomials. Let 
   $w_{mi},\ 1\leq i\leq n+1$, be sequences of points in $\C$, converging to distinct points $w_1,\ \cdots ,w_{n+1}$, respectively. 
   Suppose that for each $i,\ 1\leq i\leq n+1$, the sequence $p_m(w_{mi})$ converges. 
   Then there are $\alpha_0,\cdots ,
   \alpha _n$ such that $a_{mi}$\ converges to $\alpha_i$ for each $i,\ 0\leq i\leq n$, 
   and therefore, 
   the sequence $p_m$ converges, uniformly on compacta, to the polynomial $p(w)=
   \alpha_0+\alpha_1w+\cdots +\alpha_nw^n$.  In particular, $p_m(w_{mi})$ converge 
   to $p(w_i), 1\leq i\leq n+1$, for any sequences $w_{mi}$ converging 
   to $w_i,\ 1\leq i\leq n+1$. \rm
   \vskip 2mm
   \noindent\bf Proof.\ \rm Since $w_1,\cdots ,w_{n+1}$ are distinct the Vandermonde matrix
   $$
   V(w_1,\cdots, w_{m+1}) =\left[
   \eqalign{&1, w_1, \ \cdots \ w_1^n\cr
   &\cdots \cr
   &1, w_{n+1} ,\cdots w_{n+1}^n \cr}\right]
   $$ 
   is nonsingular and consequently for all sufficiently large $m$ the matrices
   $$
   V(w_{m1}, \cdots ,\ w_{m,n+1}) = \left[
   \eqalign{&1,w_{m1},\ \cdots \ ,w_{m1}^n\cr
   &\cdots \cr
   &1, w_{m,n+1},\cdots\ w_{m,n+1}^n\cr}\right]
   $$ are nonsingular and as $m\rightarrow\infty$, converge
   to the nonsingular matrix $V(w_1,\cdots ,w_n)$. 
   Consequently, for $m$ 
   sufficiently large, the matrices $V(w_{m1},\cdots 
   ,w_{m,n+1})^{-1}$ are well defined and, as $m\rightarrow\infty$, they 
   converge to $V(w_1,\cdots ,w_{n+1})^{-1}$. Since
   $$
   \left[\eqalign{ &p_m(w_{m 1})\cr
   &\cdots \cr
   &p_m(w_{m,n+1})\cr }\right]  = V(w_{m 1}\cdots , w_{m n}). 
   \left[\eqalign{a_{m 0}\cr\cdots\cr a_{m n}\cr}\right]
   $$
   it follows that
   $$
   \left[\eqalign{a_{m 0}\cr\cdots\cr a_{m n}\cr}\right]= 
   V(w_{m 1}\cdots , w_{m n})^{-1}
   \left[\eqalign{&p_m(w_{m, n+ 1})\cr&\cdots\cr &p_m(w_{m,n+1})\cr}\right]
   $$
   and since both factors on the right converge it follows that the columns on the left converge. This completes the proof.
   
   \centerline{$\sim$}
   
   Now, fix $x, -\ 1<x<0$, and fix distinct points $y_1,
   \cdots y_{n+1}\in (1/x , x)$. Choose a sequence $z_n,\ \Im z_n
   >0$, converging to $x$ and observe that by the nature of $bD_{z_n}$\
   for each $i,\ 
   1\leq i\leq n+1$, there is a sequence $u_{m i}\in I_{z_m}, m\in\N$ 
   that converges to $y_i$, and a sequence $w_{m i}\in J_{z_m}$ that converges 
   to $y_i$. By (10.1) $\ F(z_m, u_{m i})$ 
   converges to $\varphi (x, y_i)$ and by (10.2),\ $F(z_m, w_{m i})$ converges to 
   $\psi (x,y_i)$. Now apply Lemma 10.1 to the sequence of polynomials
   $$
   w\mapsto p_m(w)= F(z_m ,w) = a_0(z_m)+a_1(z_m)w+\cdots + a_n(z_m)w^n
   $$ 
   to see that there are 
   numbers $\alpha_i = \lim_{m\rightarrow\infty}
   a_i (z_m),\ 0\leq i\leq n$ and that for each  $i,\ 1\leq i\leq n+1$,
   $$
   \varphi (x, y_i)=\alpha_0+\alpha_1y_i+\cdots +\alpha_ny_i^n = \psi (x, y_i).
   $$
   In particular, since $y_j$ were arbitrary, it follows that
   $$
   \varphi (x,y)\equiv \psi (x,y)\ \ (1/x<y<x).
   $$
   Moreover, keeping $y_1,\cdots, y_{n+1}$ fixed we see that 
   $\alpha_0, \alpha_1,\cdots ,\alpha_n$ do not 
   depend on $z_n$ converging to $x$. This implies
   that each function $z\mapsto a_i (z)$ extends 
   continuously to $\D_+\cup\{ x\}$ and thus 
   $$
   \lim_{z\rightarrow x, \Im z>0} F(z,y)= 
   \alpha_0+\alpha_1y+\cdots \alpha_n y^n=\varphi (x,y) = \psi (x,y) .
   $$
   In the same way we get 
   $$
      \lim_{z\rightarrow x, \Im z<0} F(z,y)= 
      \beta_0+\beta_1y+\cdots \beta_n y^n=\varphi (x,y) = \psi (x,y) ,
   $$
   which, finally, implies that each $a_j,\ 0\leq j\leq n$, 
   extends holomorphically across $(-1,0)$. In the same way we get that each  
   $a_j,\ 0\leq j\leq n,$ extends holomorphically across $(t,1)$,
   \ $F$ is well defined and
   $$
   F(z,w)=a_0(z)+a_1(z)w+\cdots + a_n(z)w^n\ \ (z\in \D\setminus \{ 0,t\} ).
   $$
   \vskip 4mm
   \bf 11. Removable singularities of the coefficients at $0$ and at $t$\rm
   \vskip 2mm
   In the preceding section we proved that 
   $$
   f(z)= F(z,\overline z)= a_0(z)+a_1(z)\overline z+\cdots +a_n(z)\overline z^n
   \ \ \ (z\in\D\setminus\{ 0,t\} )
   $$
   where the functions $a_j,\ 0\leq j\leq n$, are holomorphic 
   on $\D\setminus \{ 0,t\}$. Recall that our $f$ has the property that
   $$
   \left.
   \eqalign{ &\hbox{for each\ } 
      \Gamma\in \cC_0\cup \cC_t \hbox{\ the function\ } 
      z\mapsto (z-c(\Gamma ))^nf(z)\cr
   &\hbox{extends holomorphically from \ } \Gamma .\cr }\right\}\eqno (11.1)
   $$
   \vskip 2mm
   \noindent\bf PROPOSITION 11.1\ \it Suppose that $f$ is as 
   above and suppose that (11.1) holds. Then the functions $a_j,\ 
   0\leq j\leq n$, extend holomorphically through $\D $.\rm
   \vskip 2mm
   \noindent\bf Proof.\ \rm If $\G =\bD (0,R)$ then
   $\overline z = R^2/z$ so the function
  $$
   z^n \left[ a_0(z)+{{a_1(z)R^2}\over z}+\cdots +{{a_n(z)R^{2n}}\over {z^n}}\right]
   \eqno (11.2)
   $$
   provides the holomorphic extension of $z^nf(z)$ from $\G $ to $\D 
   (0,R)\setminus \{ 0\}$ if $0<R<t$ and to $\D (0,R)\setminus\{ 0,t\} $ if $t<R<1$. 
   
   Consider first the second case. In particular, 
   when $t<R<1$ the function (11.2) has no singularity at $t$ which is the same as to say that 
   the function
   $$
    z\mapsto a_0(z)+{ {a_1(z)}\over z }{R^2}+\cdots +{ {a_n(z)}\over{z^n} }R^{2n}
   $$
   has no singularity at $t$. Let $p_j(z)$ be the singular part of 
   ${{a_j(z)}\over{z^j}}$ in the Laurent expansion around $t$. We must have
   $$
   p_0(z) + R^2p_1(z)+\cdots + R^{2n}p_n(z)\equiv 0\ \ \ (t<R<1)
   $$ 
   for $z$ in a neighbourhood of $t$ which implies 
   that $p_j(z)\equiv 0\ (0\leq j\leq n)$ so $z\mapsto a_j(z)$ 
   is holomorphic in a neighbourhood of $t$, 
   which, since $t>0$, implies that each $a_j, \ 0\leq j\leq n$, 
   is holomorphic in a neighbourhood of $t$. 
   This shows that each $a_j,\  0\leq j\leq n$, is holomorphic on
   $\D\setminus\{ 0\}$. 
   
   Now observe that for each $R,\ 0<R<1$, the function
   $$
   \eqalign{ 
   z\mapsto  z^n\bigl[ a_0(z)+a_1(z){{R^2}\over z}+\cdots +a_n(z){{R^{2n}}\over{z^n}}\bigr ] = 
   \cr
   = a_0(z)z^n+a_1(z)R^2z^{n-1}+\cdots +a_n(z)R^{2n}
   \cr}
   $$
   provides the holomorphic extension of $z^nf$ from $\bD (0,R)$ to $\D (0,R)\setminus\{ 0\}$, By our assumption this extension 
   must be holomorphic at $0$. In the Laurent series at the origin, let
   $$
   \eqalign{ &q_0\hbox{\ \ be the singular part of \ }z^na_0 \cr
   &q_1\hbox{\ \ be the singular part of \ }z^{n-1} a_1 \cr
  &\cdots \cr
   &q_n\hbox{\ \ be the singular part of \ }a_n .\cr}
   $$
   We must have 
   $$
   q_0+R^2q_1+\cdots +R^{2n}q_n \equiv 0 \ \ (0<R<1)
   $$
   which implies that
   $$
   q_0\equiv q_1\equiv\cdots\equiv q_n \equiv 0 
   $$
   and it follows that there are holomorphic functions $h_0, h_1,\cdots, h_n$ on $\D $ such that on $\D\setminus 
   \{ 0\}$ our function $f$ has the form 
   $$
   f(z)= {{h_0(z)}\over{z^n}} +{{h_1(z)}\over{z^{n-1}}}\overline z 
   +\cdots + h_n(z) \overline z ^n .
   $$
   We must now show that each $a_j$ has a removable singularity at the origin. 
   To do this, we use the fact that from 
   each circle $b\D (\lambda, \rho(\lambda ))\in C_t$ the function $f$ extends 
   holomorphically through $\D (\lambda, \rho(\lambda ))\setminus\{ \lambda \}$ 
   with a pole at $\lambda $. On 
   $\bD (\lambda, \rho(\lambda))$ we have
   $$
   \overline z = \lambda +{{\rho(\lambda )^2}\over {z-\lambda}}\ \ (0\leq \lambda<t) ,
   $$
    so our requirement will be that 
   $$
   {{h_0(z)}\over{z^n}}+ 
   {{h_1(z)}\over {z^{n-1}}}\biggl[\lambda+
   {{\rho(\lambda )^2}\over{z-\lambda}}\biggr] +\cdots + h_n(z)\biggl[\lambda +
   {{\rho(\lambda)^2}\over {z-\lambda}}\biggr]^n 
   $$
   must have no pole at the origin. This is the same as to say that
   $$
   h_0(z) +h_1(z)\biggl[ z\biggl(\lambda +
   {{\rho(\lambda^2}\over{z-\lambda}}\biggr)\biggr] + \cdots +
   h_n(z)\biggl[ z\biggl(\lambda +
   {{\rho(\lambda^2}\over{z-\lambda}}\biggr)\biggr]^n
   $$
   has a zero of order $\geq n$ at $0$ for each $\lambda,\ 0\leq \lambda <t$.
   Requiring this we should then conclude that
   $$\eqalign{
   &h_0\hbox{\ has a zero of order\ } \geq n\hbox{\ at the origin\ } \cr
   &h_1\hbox{\ has a zero of order\ } \geq n-1 \hbox{\ at the origin\ } \cr
   &\cdots \cr
   &h_{n-1}\hbox{\ has a zero of order\ } \geq 1 \hbox{\ at the origin\ } \cr}
   $$
   which is the same as to conclude that each of the functions 
   $H_0 =h_0,\  H_1=zh_1 ,\ \cdots , H_n = z^nh_n$\ has 
   a zero of order $\geq n$ at the origin. Thus, denoting
   $$\varphi (z,\lambda)= \lambda + {{\rho(\lambda )^2}\over{(z-\lambda)}}$$ we have to show that if
   $$
   H_0(z)+H_1(z)\varphi (z,\lambda)+\cdots + H_n(z)\varphi (z,\lambda)^n
   \eqno (11.3)
   $$
   has, for each $\lambda$, zero of order $\geq n$ at the origin, then the same hods for 
   each $H_j$. Suppose that (11.3) has a zero at the origin for each $\lambda $. This implies that 
   $$
   H_0(0)+H_1(0)\varphi (0,\lambda) + H_2(0)\varphi 
   (0,\lambda)^2+\cdots+H_n(0)\varphi (0,\lambda)^n\equiv 0 .
   \eqno (11.4)
   $$
   We now use the fact that
   $$
   \biggl\{\lambda - {{\rho(\lambda)^2}\over \lambda}\colon\ 0<\lambda<t\biggr\} 
   $$ 
   is a large set. This is so by the continuity: when 
   $\lambda\rightarrow 0,\ \lambda - \rho(\lambda )^2/\lambda$ is very large 
   and when $\lambda $ is near $t$,\  $\rho(\lambda ) $ is near $0$ so
   $\lambda-\rho (\lambda)^2/\lambda $ is near $t$.. Thus, (11.4) 
   implies that $H_0(0)=\cdots = H_n(0)=0$ so $ H_0(z) =zG_0(z),\cdots ,H_n(z)=zG_n(z)$ with 
   $G_i$ holomorphic. If (11.3), which now becomes 
   $$
   z\bigl[G_0(z) +G_1(z)\varphi (z,\lambda)+\cdots + G_n(z)\varphi (z,\lambda)^n\bigr] ,
   $$
   has a zero of order $\geq n$ at the origin then
   $$
   G_0(z) +G_1(z)\varphi (z,\lambda)+\cdots + G_n(z)\varphi (z,\lambda)^n
   $$ has a zero of order $\geq n-1$ at the origin. In particular,
   $$
   G_0(0)+G_1(0)\varphi (0,\lambda)+\cdots G_n(0)\varphi (0,\lambda )^n\equiv 0
   $$
   and so, after $n$ steps we deduce that each $H_j$ has zero of order $\geq n$ at the origin.
   This completes the proof of Theorem 4.1 and thus proves Theorem 1.1 in the 
   case when $\beta =0$ and $\alpha \in\D\setminus \{ 0\}$,assuming that 
   we have proved (8.1), (8.2) and (8.3).
   \vskip 4mm
   \bf 12. Proving (8.1), (8.2) and (8.3), Part 1\rm
      \vskip 2mm
      We first look at (8.1), that is, at the function $F$ on $\Omega_+$. Write
      $$
      M_+ =\{ (z,w)\in M\colon\ \Im z\geq 0\}.
      $$
      Clearly $M_+$ is a part of $b\Omega_+$. We have
      $$
      M_+ = \left[ \bigcup_{ \ z\in\DD,\ \Im z>0} \{ z\} \times bD_z 
      \right ]\bigcup\left[ \bigcup_ {z\in (0,t)}  \{ z\}\times \R  \right]  .
      $$
      We know that the function $F$ is well defined and continuous 
      on $M_+$. We also know 
      that the function $F$ is continuous on  
      $$
      \bigcup_{z\in\DD,\ \Im z>0}\overline D_z = 
      \Omega_+\bigcup \left[\bigcup _{z\in\DD, \Im z>0}\{ z\} \times bD_z \right]
      $$
      and that the extension is holomorphic on $\Omega_+$.
      We now want to show that $F$ extends continuously to 
      $\Omega_+\cup M_+\cup \Theta_+ $ where 
      $$\Theta_+ = \bigcup_{0<\eta<t}\{ \eta \} \times \{\Im \z \leq 0\}  .$$ The set $\Theta_+$ is a 
      disjoint union of halfplanes
      attached to $M_+$ along $\bigcup_ { z\in (0,t) }\{ z\}\times \R $ and is contained in $b\Omega_+$.
   In other words, we want to show that our function, holomorphic on $\Omega_+$, extends continuously to 
   $$
   \left[\bigcup _{z\in\DD,\ \Im z >0} \{z \}\times 
   \overline{D_z}\right]\bigcup\left[ \bigcup_{0<\eta<t}\{ \eta\}
   \times \{\Im\z\leq 0\}\right] 
   $$
     which is (8.1) that we want to prove.  This will imply that for each $\eta\in (0,t)$ the function 
      $\z\mapsto F(\eta,\z )$ is continuous on $\{\Im\z\leq0\}$ and holomorphic 
      on $\{\Im \z <0\}$. We will also show that it satisfies an estimate of the form 
      $$
      |F(\eta ,\z )|\leq M_\eta (1+|\z |)^n\ \ (\Im \z \leq 0)
      \eqno (12.1)
      $$
      which will, after proving the inequality also for $\Im\z\geq 0$, give (8.3). 
      To get an estimate of the form (12.1) we first prove the following
      \vskip 2mm
      \noindent \bf LEMMA 12.1\it\ Given $\eta,\ 0<\eta<t$ there are 
      an open disc $U$ centered at $\eta$ and a constant 
   $M=M(U)$ such that
   $$
   \left| {1\over{(w-i)^n}}F(z,w)\right|\leq M\ \ (z\in U,\ \Im z>0,\ w
   \in bD_z) . \eqno (12.2)
   $$
   \vskip 2mm
   \noindent\bf Remark\ \rm Note that for $z\in \D,\ \Im z>0$,\ \ $\overline {D_z}$ is contained in 
     $ \{ \Im w <0\}$ where the function $w\mapsto 1/(w-i)$ is holomorphic. 
     \vskip 2mm
     \noindent \bf Remark\ \rm By the maximum priciple the estimate (12.2) implies that 
     $$
     \left| {1\over{(w-i)^n}}F(z,w)\right|\leq M\ \ (z\in U,\ \Im z>0,\ w
     \in \overline{D_z} 
   $$
   which gives
   $$
   |F(z,w)\leq M(1+|w|)^n\ \ (z\in U,\ \Im z>0, w\in\overline{D_z})
   $$
   and which, after proving  the desired continuous extendibility, gives (12.1).
   \vskip 2mm
   \noindent \bf Proof.\ \rm Note first that 
   $w\mapsto 1/(w-i)$ is holomorphic on $\{ \Im w\leq 0\} $ and 
   satisfies $|1/(w-i)|\leq 1\ (\Im w\leq 0)$. Fix $\eta, \ 0<\eta<t$. There are an open 
   disc $U$ centered at $\eta$ and $T_0$ and $R_0$ such that 
   $$
   U\cap\DD(T,\rho (T)) =\emptyset\ \ (T_0<T<\eta)\ \ \hbox{and}\ \  
   U\cap \DD (0,R)=\emptyset \ \ (0<R<R_0) .
   $$
   In other words, if $z\in U$ and $\{ z\}\times \C$ meets $
   \Lambda (T,\rho (T))\cup b\Lambda (T,\rho (T))  $
   then necessarily $0\leq T\leq T_0$, and if $\{ z\}\times \C$ meets 
   $\Lambda (0,R)\cup b\Lambda (0,R)$ then necessarily $R_0\leq R\leq 1$.
   
   Recall that 
   $\Lambda (T,\rho (T))=\{ (\z, T+\rho (T)^2/(\z -T))\colon\ 0<|\z - T|<\rho (T)\} $
   and $\Lambda (0,R)=\{ (\z , R^2/\z )\colon\ 0<|\z |<R\}$ and notice that for 
   $w=T+\rho (T)^2/(z-T)$ we have $1/(w-i)= (z-T)/\left[ (T-i)(z-T)+\rho(T)^2\right]$ and 
   for $w=R^2/z$ we have $1/(w-i)= z/(R^2-iz)$.
   
   Choose $\varepsilon >0$ so small that $\varepsilon < \rho (T_0)^2/(
   2\sqrt 2)$ and $\varepsilon <R_0^2/2$. Assume that $z\in U,\ \Im z>0$. 
   If $\{ z\}\times\C$ meets $\Lambda(T,\rho (T))\cup b\Lambda(T,\rho (T))$ 
   then 
   $0\leq T\leq T_0$ and $\rho (T)\geq\rho (T_0)$ so if $|z-T|<\varepsilon $ then 
   $$
   |\rho (T)^2+(T-i)(z-T)|\geq 
   \rho (T_0)^2-\sqrt 2\varepsilon \geq \rho (T_0)^2/2 .
   $$
   If $\{ z\} \times\C $ meets $\Lambda (0,R)\cup b\Lambda (0,R)$ it 
   follows that $1\geq R\geq R_0$ and consequently, 
   for $|z|<\varepsilon$ we have 
   $|R^2-iz|\geq R_0^2-\varepsilon \geq R_0^2/2.$
   
   Recall that on $\Lambda (T,\rho (T)$ we have
   $$
   F\left(z,T+{{\rho (t)^2}\over{z-T}}\right) = 
   {{d_{-n}(T)}\over {(z-T)^n}}\cdots +{{d_{-1}(T)}\over{z-T}} + h_T\left (
       {{z-T}\over {\rho(T)}}\right ) 
       $$
       where, by Lemma 4.1,  the functions $d_{-1},\cdots ,d_{-n}$ are continuous 
       and $h_T$ are functions from the disc algebra, continuously depending 
       on $T$. On $\Lambda (0,R)$ we have
   $$ F\left( z, {{R^2}\over z}\right) =
       {{c_{-n}(R)}\over
       {z^n}}+\cdots + {c_{-1}(R)\over z}+g_R\left( z\over R\right)
       $$
   where, by Lemma 4.1, the functions $c_{-1},\cdots ,c_{-n}$ are continuous and $g_R$ are functions from 
   the disc algebra, continuously depending on $R$.  
   
   Fix $z\in U,\ \Im z >0$, and fix $T$ such that $\{ z\}\times\C$ meets
   $\Lambda (T,\rho (T)\cup b\Lambda (T,\rho (T))$. 
   Let $(z,w(z)) \in \Lambda (T,\rho (T)\cup b\Lambda (T,\rho (T))$.
   We know that this implies that $0\leq T\leq T_0$ so $1\geq \rho (T)\geq
   \rho (T_0)$. Suppose first that $|z-T|<\varepsilon$. Then
   $$
   \eqalign{
   &\left|
   \left( {1\over{w(z)-i}} \right)^n  F(z,w(z)) \right| = \cr
   &= \left| \left[ {{z-T}\over{(T-i)(z-T)+\rho(T)^2}} 
   \right]^n.\left[ {{d_{-n}(T)}\over 
   {(z-T)^n}}\cdots +{{d_{-1}(T)}\over{z-T}} + h_T\left (
       {{z-T}\over {\rho(T)}}\right )\right]\right| = \cr
    &=\left| {1\over{(T-i)(z-T)+\rho(T)^2}} \right|^n.\cr  
    &. \left| d_{-n}(T) + d_{-n+1}(T)(z-T)
   +\cdots d_{-1}(T)(z-T)^{n-1}+ 
   (z-T)^n .
   h_T\left (
       {{z-T}\over {\rho(T)}}\right)\right| \leq\cr
   &\leq \left[{2\over{\rho(T_0)^2}}\right]^n.\left[ |d_n(T)|+
   \cdots +|d_{-1}(T)|\varepsilon^{n-1} + \varepsilon^n \left| h_T
   \left( {{z-T}\over {\rho(T)}}\right)\right| \right]\cr}
   $$
   Since $|1/(w(z)-i)|\leq 1$ it follows that if $|z-T|\geq \varepsilon$ then 
   $$
   \left|
   \left( {1\over{w(z)-i}} \right)^n  F(z,w(z)) 
   \right| \leq {{|d_{-n}(T)|}\over{\varepsilon^n}}+\cdots +{{|d_{-1}(T)|}\over\varepsilon} +
   \left| h_T
   \left( {{z-T}\over {\rho(T)}}\right)\right| .
   $$
   Since the functions $ d_j,\ -n\leq j\leq -1$ are continuous and since $ h_T$, the 
   functions in the disc algebra, continuously depend on $T$ the right hand sides of both inequalities 
   above are uniformly bounded for $0\leq T\leq T_0$. It follows that the function $(z,w)\mapsto F(z,w)/(w-i)^n$
   is bounded on $\bigcup_{z\in U,\ \Im z>0}\{ z\} \times \lambda_z$.
   
   Now, let $z\in U,\ \Im z >0$ and fix $R$ such that $\{ z\}\times\C$ 
   meets $\Lambda (0,R)\cup b\Lambda (0,R)$. We 
   know that this implies that $R\geq R_0$. Let 
   $( z,w(z)\in\Lambda (0,R)\cup b\Lambda (0,R)$. As above, suppose
   first that $|z|<\varepsilon$. Then
   $$
   \eqalign{
   & \left| \left( {1\over{w(z)-i}}\right)^nF(z,w(z))\right| = \cr
   &= \left|\left( {z\over{R^2-iz}}\right) ^n\left[ {{c_n(R)}\over{z^n}}
       + \cdots {{c_{-1}(R)}\over z} +
       g_R\left ( {z\over R}\right )\right ]\right| = \cr
    &   =\left|{1\over {(R^2-iz)^n}}\left[ c_{-n}(R) +\cdots +c_{-1}(R)z^{n-1} 
       +z^n g_R\left( {z\over R} \right )\right ] \right| \leq \cr
    & \leq \left( {2\over{R_0^2}}\right)^n.\left[ |c_{-n}(R)|+\cdots +|c_{-1}(R)| + 
    \left| g_R\left( {z\over R}\right)\right|\right]\cr     }
       $$
     Since $|1/(w(z)-i)|\leq 1$ it follows that if $|z-T|\geq \varepsilon$  then,  as above,
     $$
     \left| \left( {1\over{w(z)-i}} \right)^nF(z,w(z)) \right| \leq 
     \left [ {{|c_{-n}(R)|}\over{\varepsilon ^n}}+\cdots +{{|c_{-1}(R)|}\over\varepsilon} + 
     \left| g_R\left ( {z\over R}\right)\right|\right] . 
     $$
   Since the functions $ c_j,\ -n\leq j\leq -1$ are continuous and since $ g_R$, the 
   functions in the disc algebra, continuously depend on $R$ the right hand sides of both inequalities 
   above are uniformly bounded for $R_0\leq R\leq 1$. It follows that the function 
   $(z,w)\mapsto F(z,w)/(w-i)^n$
   is bounded on $\bigcup_{z\in U,\ \Im z>0}\{ z\} \times I_z$ which completes the proof of Lemma 12.1. 
   \vskip 4mm
    13.\ \bf Completing the proof of (8.1), (8.2) and (8.3) \rm
   \vskip 2mm
   We use the map $w\mapsto 
      \varphi (w) = 1/(w-i)$  
      which maps $\{ \Im\z\leq 0\}$ homeomorphically 
      to $\overline \P\setminus\{ i\} $ where
       $\P$ is the disc in $\C$ centered at $i/2$ of 
       radius $1/2$ and which maps $\{ \Im w <0\}$ biholomorphically to $\P$. 
       Note that $\varphi ^{-1}(W)=1/W +i$. Thus, the map 
        $(z,w)\mapsto \Phi (z,w)= (z,\varphi (w))$ 
       maps $\C \times\{ \Im w\leq 0\}$ homeomorphically to 
       $\C\times (\overline \P\setminus\{ 0\})$ 
       and $\C\times \{ \Im w <0\}$ biholomorphically to $\C\times\P$.
       
       The map $\Phi$ maps the domain $\Omega_+$ biholomorphically to the domain
       $\Sigma_+ = \Phi (\Omega_+)= \bigcup_{z\in\D,\ \Im z>0}\{ z\}\times E_z$
       where for each $z\in\D,\ \Im z>0,\ E_z =\varphi (D_z) $ is a domain in $\P$ whose 
       boundary consists of two circular arcs with endpoints $\varphi (\overline z)$ and $\varphi (1/z)$. When $z,\Im z>0$, converges
       to a point $\eta\in (0,t)$, the endpoints converge to the
       points $\varphi (1/\eta)$ and $\varphi (\eta)$ on $b\P$ and the 
       domains $E_z$ converge to $\P$. Clearly 
       $$
       \Phi (M_+)=\left[\bigcup_{z\in\DD,\Im z>0}\{z\}\times 
       bE_z\right]\bigcup\left[\bigcup_{z\in (0,t)}
       \{ z\} \times (b\P\setminus\{ 0\})\right]
       .
       $$
   
   We want to prove that $F$ extends continuously 
   to $[U\cap\R]\times \{ \Im w\leq 0\}$. This will be done if we prove that the function 
   $W\mapsto F(z,1/W+i)$, which is continuous on 
   $$\left[\bigcup_{z\in U,\Im z >0} \{ z\}\times bE_z\right]\bigcup \left[ 
   \bigcup_{z\in U\cap\R}\{ z\}\times [b\P\setminus\{ 0\}]\right]$$ 
   and 
   which is continuous on 
   $$\bigcup_{z \in U, \Im z>0}\{ z\}\times \overline{E_z} $$ 
   and 
   holomorphic on \ $\bigcup_{z\in U, \Im z>0}\{ z\}\times E_z$, 
   extends continuously to  
   $\bigcup_{z \in U\cap\R}\{ z\} \times [\overline\P\setminus\{ 0\}]$. 
   By what we have shown above the function 
   $$
   G(z,W)=\left\{\eqalign{
   &W^{n+1}F(z, 1/W+i)\ \ (W\not=0) \cr
   & 0\ \ (W=0)\cr}\right.
   $$
   is continuous on 
   $$
   \left[\bigcup_{z\in U,\Im z >0} \{ z\}\times bE_z\right]\bigcup \left[ 
   \bigcup_{z\in U\cap\R}\{ z\}\times b\P\right]
   $$ and continuous on 
    $\bigcup_{z \in U, \Im z>0}\{ z\}\times \overline{E_z} $ and holomorphic on 
   $\bigcup_{z\in U, \Im z>0}\{ z\}\times E_z$. 
      The continuous extendibility of $F$ which we want to prove will be proved once we have 
      shown that our function $G$ extends continuously to
      $$
      \left[\bigcup_{z\in U,\ \Im z>0}\{ z\} \times\overline{E_z}\right]\bigcup 
      \left[\bigcup_{\eta\in U\cap\R}\{\eta\}\times\overline \P \right].
      $$
      With no loss of generality assume that $0\in E_z$ for all $z\in U,\ \Im z > 0$. We understand that 
      $E_z=\P$ for $z\in U\cap\R$. 
      The domains $E_z$ change continously with $z$ in the sense that $bE_z$ change continuously with $z$ in
      Frechet's sense [Ts, p.383] which, by a theorem of Courant [Ts, p.383] implies that if for each 
      $z\in U,\ \Im z\geq 0$,\ $\psi_z$ is the conformal map from $\D$ to $E_z$,\ $\psi_z(0)=0,\ \psi_z^\prime (0) >0$, 
      these maps change continuously with $z$, uniformly on $\DD$. It follows that if we set 
      $\Psi (z,w) = (z,\psi_z(w))$ then $\Psi $ maps $\{ z\in U,\ \Im z\geq 0\}\times \DD$ 
      homeomorphically onto $\bigcup_{z\in U,\ \Im z\geq 0}\{ z\}\times\overline{E_z}$. 
      Recall that the function $G$ is continuous on 
      $\bigcup _{z\in U\ \Im z\geq 0}\{ z\}\times bE_z$ and is also continuous on 
      $\bigcup_{z\in U,\ \Im z>0}\{ z\}\times \overline{E_z}$ and holomorphic on 
      $\bigcup _{z\in U,\ \Im z>0}\{ z\}\times E_z$.
      
      Let $J=G\circ\Psi$. The function $J$ is continuous on 
      $\{ z\in U,\Im z >0\}\times\DD$ and holomorphic on each $\{ z\}\times\D,\ z\in 
      U,\ \Im z>0$, and is also continuous on $\{ z\in U , \ \Im z\geq 0\}\times\bD$. 
      By the continuity it follows that for each $z\in (0,t)$ the function $J$ extends holomorphically through 
      $\{ z\}\times\D$ and if we define the extension $\tilde J$ so that on each 
      $\{ z\}\times \D$ it is the holomorphic extension of $J$ from $\{ z\}\times \bD $,
      then so extended function will be, by the continuity of $J$ on
      $\{ z\in U,\ \Im z\geq 0\} \times\bD$ and by the maximum principle, continuous on 
      $\bigcup_{z\in U,\ \Im z\geq 0}\{ z \}\times \DD $. It 
      follows that $\tilde J\circ\Psi^{-1}$ provides the necessary extension of 
      $G$ to $\bigcup_{ z\in U,\ \Im z\geq 0}\{ z\}\times \overline {E_z}$. 
      
      This completes the proof of (8.1), and "half" of (8.3). In the same way we prove 
      (8.2) and the "other 
      half" of (8.3). Theorem 4.1 is proved. This proves Theorem 1.1 in the case (2.2). 
      
      The proof of Theorem 1.1 in the cases (2.3) and (2.4) is almost the same. 
      The case (2.3) is the
      limiting case of (2.2) when $t$ tends to $1$.  For $z\in\D,\ z\not\in \R$, 
      the domain $D_z$ is now 
      bounded by $I_z$, the segment joining $\overline z$ and $1/z$, and 
      by $\lambda_z$, the arc of the circle passing thhrough $\overline z,\ 1/z$ and $1$, with
      endpoints $\overline z$ and $1/z$, which does not contain $1$. 
      Folding occurs only on the interval $(-1,0)$ and there is only the singularity at the origin to remove in 
      proving the analyticity of the coefficients $a_0,\cdots ,a_n$ on $\Delta$. 
      The case (2.4) is even simpler. 
      The domain $D_z$ is now bounded by $\lambda_z$ and $\mu_z$ where $\mu_z$ 
      is the arc on the circle passing through $\overline z,\ 1/z$ and 
      $-1$, with endpoints $\overline z$ and $1/z$, which does not contain $-1$.
      There is no folding and no singularity to remove. 
   Theorem 1.1 is proved.
   \vskip 4mm
   \bf 14.\ Remarks and open questions
   \vskip 2mm
   \noindent \bf Remark \ \rm In Theorem 1.1 the continuity of $f$ at the boundary is
   essential as shown by the example 
   $$
   f(z) ={{z(z-1/2)}\over{1-|z|^2}}\ \ (z\in\D)
   $$
   of a function which extends holomorphically from every circle $\Gamma\in \cC_0\cup \cC_{1/2}$ yet 
   it is not holomorphic on $\D$.
   \vskip 2mm
   \noindent
   \bf Remark \ \rm Note that if $f$ is continuous on $\overline\D$ and of the form (1.1) then
   $a_j,\ 0\leq j\leq n,$ need not be continuous on $\overline\D$. To see this, let $g$ be a bounded 
   holomorphic function on $\D$ that does not extend continuously to $\overline\D$. Then 
   $f(z)= (1-|z|^2)g(z) = g(z) - [zg(z)]\overline z$ extends continuously 
   through $\overline\D$ yet the functions $g$
   and $zg$ do not. 
   \vskip 2mm
   \noindent \bf Example 14.1\  \rm The function $g(z,w)= |w|^2\ ((z,w)\in b\B )$
   extends holomorphically into $\B$ along each complex line
   passing through the origin and along each complex line parallel to one 
   of the coordinate axes that meets $\B$. Let $M$ be 
   the Moebius transform in $\C^2$ which maps the origin to the point 
   $a=(1/2, 1/2)$ [Ru]. $M$ maps the the complex lines through the origin 
   to $\L (a)$, the complex lines parallel to the $z-$axis to $\L (b)$ 
   and the complex lines parallel to the $w-$axis to $\L (c)$ where 
   $b, c$ are contained in $\C^2\setminus\overline\B$ and satisfy
   $<a|b>=<a|c>=1,\ <b|c>\not= 1.$  The function $f=g\circ M^{-1}$ is real analytic on 
   $b\B$ and extends into $\B$ along each complex line
   $L\in\L (a)\cup\L (b)\cup\L (c)$ yet it does not 
   extend holomorphically into $B$, and so, by the 
   main result of [G3], \ $a$ is the only point in $\B$ 
   such that $f$ extends holomorphically into
   $\B$ along each complex line through this point. 
   \vskip 2mm
   \noindent This example shows that in Theorem 1.2 and Corollary 1.3 one cannot drop the assumption that $<a|b>\not=1$ if one of the points $a, b$ is in $\B$ and the other in $\C^2\setminus\overline\B$. It 
   is easy to see that the points $a, b, c$ do not lie on the same complex line so the example shows also that there are triples $a, b, c$ of points in $\C^2$, 
   not lying on the same complex line, such that $\L (a)\cup\L (b)\cup \L (c)$  is 
   not a test family for holomorphic extendibility for $\cC^\infty (b\B)$,
   \vskip 2mm
   \noindent\bf Example 14.2\ \rm Let $$
g(z,w)=\left\{\eqalign{&{z^{k+2}\over\overline z}\ \ \ (z\not=0)\cr
&0\ \ \ (z=0)\cr}\right.
$$ This is a function of class $\cC^k$ on $b\B$ which extends holomorphically 
into $\B$ along any complex line which meets $\{0\}\times\D$. In particular,
it extends holomorphically into $\B$ along every
complex line that is parallel to the $z-$axis. Let $M$ be 
the Moebius map in $\C^2$ that maps the origin to the point $(1/2, 0)$. 
The function $f=g\circ M^{-1}$ extends holomorphically into $\B$ along 
any complex line that meets $L\cap \B$ where $L=\{ (z,w)\colon \ z=1/2\}$ and
along each complex line passing through $(2,0)$ yet it does not extend 
holomorphically through $\B$. Choosing $a, b\in L\cap \B$ and $c=(2,0)$ 
we have $<a|b>=<a|c>=1$ which shows that in Corollary 1.3 we cannot drop
the assumption that one of the numbers 
$<a|b>,\ <a|c>$ is different from $0$.
\vskip 2mm
If $E$ is a complex line that misses $\overline\B$ then
   given a finite set $\{ a_1,a_2,\cdots ,a_n\} \subset E$ there
   is a real analytic function on $b\B$ which extends holomorphically  
   into $\B$ along every $L\in \L (a_1)\cup\L (a_2)\cup\cdots\cup\L (a_n)$, yet it 
   does not extend holomorphically through $\B$ [KM, G3]. The situation in 
   the case when $E$ is tangent to $b\B$ is unclear. If 
   $E$ meets $\B$ the situation is different for functions of class
   $\cC^\infty $: the lines through two points $a_1,\ a_2$ suffice 
   provided that 
   $<a_1|a_2>\not= 1$ in the case that one of the points $a_1, a_2$ is in $\B$ and the  other in 
   $\C^2\setminus\overline\B$ [G3]. The situation is different for functions of class $\cC^k$:
   Given numbers $t_1, t_2,\cdots , t_n,\ 1<t_1<t_2<\cdots <t_n<\infty$ the function
   $$
   (z,w)\mapsto \left\{ 
   \eqalign {
   &{{w^{k+2}}\over{\overline w}}\left( z-{1\over{t_1}}\right)\cdots\left( z-{1\over{t_n}}\right)\ \ \ (w\not= 0) \cr
   &0\ \ \ (w=0)\cr}\right.
   $$
   is a function of class $\cC^k$ on $b\B$ which extends holomorphically into  $B$ 
   along each complex line which meets $\D\times\{ 0\} $ as well as along each complex line $L\in \L ((t_1,0))\cup 
   \L ((t_2, 0))\cup\cdots\cup \L ((t_n, 0)) $ yet does not extend 
   holomorphically through $\B$.
   
   We conclude with the following open
   \vskip 2mm
   \noindent \bf QUESTION\ \it Suppose that $a, b, c\in \C^2$ do not lie on the 
   same complex line and assume that $\Lambda (a,b)\cap\overline\B = 
   \Lambda (b,c)\cap\overline\B = \Lambda (a,c)\cap\overline\B = \emptyset$. 
   Is $\L (a)\cup \L (b)\cup \L (c)$ 
   a test family for holomorphic extendibility for $\cC (b\B )$? 
   If not, is this family a test family for
   holomorphic extendibility for $\cC^\infty (b\overline\B )$?\rm

  \vskip 4mm
  \noindent\bf Acknowledgement\rm\ 
  The author is indebted to M.\ Agranovsky for providing the example in Section 3.\ \ A part of this paper was
  written in September 2010 during the author's stay at the
  Institute of Mathematics, University of Bern, Switzerland.  He wishes to thank 
  Frank Kutzschebauch for the invitation.
  \vskip 3mm
  This paper was supported in part by the ministry of Higher Education, Science and Technology of Slovenia through the research
  program Analysis and Geometry, Contract No. P1-02091 (B).
 \vfill
 \eject

 \centerline{\bf REFERENCES}
 \vskip 5mm
 \noindent [A1]\ M.\ L.\ Agranovsky:\ Holomorphic extension from the unit
 sphere in $C^n $ into complex lines passing through a finite set.
 
 \noindent To appear in J.\ d'Analyse Math.\ \   http://arxiv.org/abs/0910.3592
 \vskip 2mm
 \noindent [A2]\ M.\ L.\ Agranovsky:\ Characterization of polyanalytic functions by meromorphic extensions into
 chains of circles.
 
 \noindent To appear in J.\ d'Analyse Math.\ \   http://arxiv.org/abs/0910.3578
 \vskip 2mm
 \noindent [A3]\ M.\ L.\ Agranovsky:\ Boundary Forelli theorem for the sphere in 
 $\C^n$ and $n+1$ bundles of complex lines.
 
 \noindent http:/arxiv.org/abs/1003.6125
 \vskip 2mm 
 \noindent [AG]\ M.\ L.\ Agranovsky and J.\ Globevnik:\ Analyticity
 on circles for rational and real-analytic functions of two real variables.
 
 \noindent J.\ d'Analyse Math.\ 91 (2003) 31-65
 \vskip 2mm
 \noindent [AV]\ M.\ L.\ Agranovsky and R.\ E.\ Val'ski:\ Maximality of invariant algebras of functions.
 
 \noindent Sibirsk.\ Mat.\ Zh.\ 12 (1971) No.\ 1, 3-12 
 \vskip 2mm
 \noindent [B1]\ L.\ Baracco: Holomorphic extension from the sphere to the ball.
 
 \noindent http://arxiv.org/abs/0911.2560
 \vskip 2mm
 \noindent [B2]\ L.\ Baracco: Separate holomorphic extension along lines 
 and holomorphic extension of a continuous function from the sphere to 
 the ball: solution of a conjecture by M. Agranovsky 
 
 \noindent http://arxiv.org/abs/1003.4705
 \vskip 2mm
 
 \noindent [G1] \ J.\ Globevnik:\ Holomorphic extensions from open families of circles.
 
 \noindent Trans.\ Amer.\ Math.\ Soc.\ 355 (2003) 1921-1931
 \vskip 2mm
 \noindent [G2]\ J.\ Globevnik:\ Analyticity of functions analytic on circles
 
 \noindent Journ.\ Math.\ Anal.\ Appl. 360 (2009) 363-368
 \vskip 2mm
 \noindent [G3]\ J.\ Globevnik:\ Small families of 
 complex lines for testing holomorphic extendibility
 
 \noindent To appear in Amer.\ J\ Math.\ \ \ http://arxiv.org/abs/0911.5088
 \vskip 2mm
 \noindent [H]\ K.\ Hoffman;\it \ Banach spaces of analytic functions \rm
 
 \noindent Englewood Cliffs (N.J.), Prentice Hall 1962
 \vskip 2mm
 \noindent [KM]\ A.\ M.\ Kytmanov, S.\ G.\ Myslivets:\ On families of complex lines sufficient for holomorphic extension.
 
 \noindent Mathematical Notes 83 (2008) 500-505 (translated from Mat.Zametki 83 (2008) 545-551
 \vskip 2mm
 \noindent [L]\ H.\ Lewy: On the local character of the
 solutions of an atypical linear differential equation in three
 variables and a related theorem for regular functions of two complex variables.
 
 \noindent Ann. of Math.\ 64 (1956) 514-522
 \vskip 2mm
 \noindent [R]\ H.\ Rossi: A generalization of a theorem of Hans Lewy.
 
 \noindent Proc.\ Amer.\ Math.\ Soc.\ 19 (1968) 436-440
 \vskip 2mm
 \noindent [Ru]\ W.\ Rudin:\ \it Function theory in the unit ball of $C^n$.\rm
 
 \noindent Springer, Berlin-Heidelberg-New York 1980
 \vskip 2mm
 \noindent [S]\ E.\ L.\ Stout: The boundary values of holomorphic functions 
 of several complex variables. 
 
 \noindent Duke Math.\ J.\ 44 (1977) 105-108
 \vskip 2mm
 \noindent [T1]\ A.\ Tumanov: A Morera type theorem in the strip.
 
 \noindent Math.\ Res.\ Lett.\ 11 (2004) 23-29
 \vskip 2mm
 \noindent [T2]\ A.\ Tumanov: Testing analyticity on circles.
 
 \noindent Amer.\ J.\ Math.\ 129 (2007) 785-790
 \vskip 2mm
 \noindent [Ts]\ M.\ Tsuji:\ \it Potential theory in modern function theory.\rm
 
 \noindent Maruzen, Tokyo 1959
 
 \vskip 8mm
 \noindent Institute of Mathematics, Physics and Mechanics
 
 \noindent University of Ljubljana, Ljubljana, Slovenia
 
 \noindent josip.globevnik@fmf.uni-lj.si

 \bye